\mag=1200
\input amstex
\documentstyle{amsppt}
\NoBlackBoxes
\TagsOnRight

\topmatter
\author
Taras E. Panov
\endauthor
\title
Calculation of Hirzebruch genera for manifolds acted on by the
group~$\Bbb Z/p$ via invariants of the action
\endtitle
\thanks
Partially
supported by the Russian Foundation for Fundamental Research,
grant no.~96-01-01414.
\endthanks
\address
Department of Mathematics and Mechanics, Moscow State
University, 119899 Moscow, Russia
\endaddress
\email
{\tt tpanov\@mech.math.msu.su}
\endemail
\abstract
We obtain general formulae expressing Hirzebruch genera of a manifold with
$\Bbb Z/p$-action in terms of invariants of this action (the sets of
weights of fixed points). As an illustration, we
consider numerous particular cases of well-known genera,
in particular, the elliptic genus. We also describe the connection with
the so-called Conner--Floyd equations for the weights of fixed points.
\endabstract
\endtopmatter
\rightheadtext{Calculation of Hirzebruch genera for manifolds}

\define\td{\operatorname{td}}
\define\rk{\operatorname{rk}}
\define\SL{\operatorname{SL}}
\define\inda{\operatorname{ind}}
\define\ida{\operatorname{id}}
\define\Ima{\operatorname{Im}}
\define\tr{\operatorname{tr}}
\define\Tr{\operatorname{Tr}}
\define\Xa{\operatorname{X}}
\define\ch{\operatorname{ch}}
\define\sh{\operatorname{sinh}}
\define\sn{\operatorname{sn}}
\define\arsh{\operatorname{arcsinh}}
\define\th{\operatorname{tanh}}
\define\arth{\operatorname{arctanh}}
\define\sign{\operatorname{sign}}
\define\under{\underline}
\define\ov{\overline}
\define\barCalJ{\,\,\,\bar{\!\!\!\Cal J}}

\document

\head
Introduction
\endhead

In this paper we obtain general formulae expressing Hirzebruch genera of a
manifold acted on by $\Bbb Z/p$ with finitely many fixed points or fixed
submanifolds with trivial normal bundle via invariants of this action. We
also describe the connection with the so-called Conner--Floyd equations for
the weights of fixed points.

Actions of $\Bbb Z/p$ were studied in~\cite{12}, \cite{13},
\cite{11},~\cite{8}, where the so-called Conner--Floyd equations were deduced
within cobordism theory (see formulae~\thetag{31},~\thetag{32}).  These
equations form necessary and sufficient conditions for a sets of elements of
$\Bbb Z/p$ to be the set of weights of some $\Bbb Z/p$-action (see~\S\,3 for
the definition).  Two approaches for the calculation of Hirzebruch genera of
a stably complex manifold with a $\Bbb Z/p$-action were proposed in~\cite{5}.

The first approach is based on the application of the Atiyah--Bott--Lefschetz
fixed point formula~\cite{1}, and so for its realization it is necessary to
have an elliptic complex of bundles that are associated to the tangent
bundle of the manifold. The Atiyah--Bott--Lefschetz formula obtained
in~\cite{1} generalizes the classical Lefschetz formula for the number of
fixed points and enables us to calculate the equivariant index of an elliptic
complex of bundles over a manifold by means of certain contribution functions
of the fixed submanifolds (see the details below). In particular, if an
operator acts on a manifold with finitely many fixed points, the
corresponding equivariant index can be expressed in terms of the fixed point
weights.  It was shown in~\cite{5} how to express the Todd genus, which is
the index of an elliptic complex (namely, the Dolbeault complex) over a
manifold with $\Bbb Z/p$-action, via the equivariant index of the same
complex for the action of the generator of $\Bbb Z/p$. This equivariant index
enters into the Atiyah--Bott--Lefschetz formula.  In this way one deduces the
formulae expressing the Todd genus in terms of the weights of fixed points
for this $\Bbb Z/p$-action. The formulae obtained by this method contain the
number-theoretical trace of a certain algebraic extension of
fields of degree $(p-1)$.  In this paper, we use the same approach to obtain
the formulae for other genera of manifolds with $\Bbb Z/p$-action: the
signature (or the $L$-genus), the Euler number, the $\hat A$-genus, the
general $\chi_y$-genus, and the elliptic genus. By this method we also obtain
some general equations (see.~\S\,5) for an arbitrary Hirzebruch genus having
the property to be the index of a certain elliptic complex of bundles
associated to the tangent bundle of a manifold with $\Bbb Z/p$-action.

In this paper we use the generalized Lefschetz formula in the somewhat
different formulation, stated in~\cite{3}. This formula and especially the
``recipe" it suggests for calculating the equivariant index of a complex via
contribution functions of fixed points (see.~\S\,5) are more convenient for
applications than the formula from~\cite{1} used in~\cite{5}. The generalized
Lefschetz formula was deduced in~\cite{3} from the cohomological form of the
Atiyah--Singer index theorem, which was also proved there. We apply both
formulae: some of the results (see~\S\,4) we obtain are based on the ``old"
Lefschetz formula of~\cite{1}, while others use the ``recipe" in~\S\,5 based
on the formula from~\cite{3}.

Another approach to the equations for Hirzebruch genera, also
proposed in~\cite{5}, is based on an application of cobordism theory in the
same way as in the derivation of the Conner--Floyd equations
in~\cite{12},~\cite{13}. In~\cite{5}, the authors give a
formula expressing the $\bmod\,p$ cobordism class of a stably complex
manifold with a $\Bbb Z/p$-action in terms of invariants of the action.
Then they show that the difference between the two
formulae for the Todd genus (obtained by the first and second methods)
is exactly the sum of the Conner--Floyd equations for the Todd genus. In this
paper we show (see Theorem~7.1) that the difference between the two seemingly
different formulae deduced by these two methods for an arbitrary
genus is a weighted sum (with integer coefficients) of the Conner--Floyd
equations for this genus.

A case of particular interest is that of the so-called elliptic genus.  In
Witten's papers, a certain invariant was assigned to each oriented
$2n$-dimensional manifold~$M^{2n}$. This invariant is the equivariant index
of the Dirac-like operator for the canonical action of circle $S^1$ on
manifold's loop space. S.~Ochanine ~\cite{14} showed that this index is a
Hirzebruch genus corresponding to the elliptic sine; which led to the term
``elliptic genus".  In~\cite{2},~\cite{9} and other papers, the rigidity
theorem for the elliptic genus of manifolds with $S^1$-action was proved.
This theorem states that if we regard the equivariant elliptic genus
$\varphi_{S^1}(M)$ of such a manifold as a character of the group~$S^1$, then
$\varphi_{S^1}(M)$ is the trivial character and is equal to the elliptic
genus $\varphi(M)$. At the same time, the elliptic genus takes its values in
the ring $\Bbb Z\left[\frac12\right][\delta,\varepsilon]$, and its value on
any manifold~$M^{2n}$ is a modular form of weight~$n$ on the subgroup
$\Gamma_0(2)\subset\SL_2(\Bbb Z)$ (cf.~\cite{7}). In this paper we obtain
formulae for the elliptic genus of a manifold with a $\Bbb Z/p$-action having
finitely many fixed points. A summary of results of this part of the paper
has already been published in~\cite{15}. As an application we deduce certain
relations between the Legendre polynomials by applying our formulae to a
special action of~$\Bbb Z/p$ on~$\Bbb CP^n$ (see~\S\,8).

In the remaining part of this article (see~\S\,9) we
generalize our constructions to the case of $\Bbb Z/p$-actions with fixed
submanifolds whose normal bundles are trivial.

\head
\S\,1. Necessary information about Hirzebruch genera
\endhead

Let $M^{2n}$ be a manifold with a complex structure in its stable tangent
bundle, that is, there is~$k$ such that $TM\oplus2k$ is a complex bundle.
We write the total Chern class of the tangent bundle $TM$ as
$$
c(TM)=1+c_1(M)+c_2(M)+\cdots+c_n(M)=(1+x_1)(1+x_2)\cdots(1+x_n).
$$
This means that the total Chern class of~$M$ is written as the product of
the Chern classes of ``virtual" line bundles whose sum gives $TM$. Therefore,
$c_i(M)$ is the~$i$th elementary symmetric function in $x_1,\dots,x_n$.

To each series of the form $Q(x)=1+\cdots$ with coefficients in a certain
ring~$\Lambda$ there corresponds the Hirzebruch genus
$\varphi_Q(M^{2n})=\bigl(\,\prod_{i=1}^nQ(x_i)\bigr)[M^{2n}]$
(see~\cite{6}). Along with $Q(x)$, we introduce $f(x)=x/Q(x)$ and
$g(u)=f^{-1}(u)$. Each Hirzebruch genus~$\varphi$ gives rise to a formal group
law $F_{\varphi}(u,v)=g_{\varphi}^{-1}\bigl(g_{\varphi}(u)+
g_{\varphi}(v)\bigr)$ with logarithm $g_{\varphi}(u)$ (see~\cite{5}).
The corresponding power system
$[u]^{\varphi}_n=g_{\varphi}^{-1}\bigl(ng_{\varphi}(u)\bigr)$
(the $n$th power in the formal group law $F_{\varphi}$).

Hirzebruch genera for real orientable manifolds~$M^{4n}$ are defined
similarly. Here we replace the Chern classes~$c_i$ by the Pontryagin
classes~$p_i$, and~$x_i$ by~$x_i^2$, that is,
$$
p(TM)=1+p_1(M)+p_2(M)+\cdots+p_n(M)=(1+x_1^2)(1+x_2^2)\cdots(1+x_n^2).
$$
Hence we now have $Q(x^2)$ instead of $Q(x)$.

\medskip

The formulae obtained in this paper refer mainly to the following
Hirzebruch genera.

\noindent 1. The universal genus~$\under{\varphi}$ corresponds to the identity
homomorphism $\ida\: \Omega_U\otimes\Bbb Q\to\Omega_U\otimes\Bbb Q$.
The corresponding formal group law of ``geometric cobordisms"
$\under F(u,v)$ (see~\cite{12}) is universal, that is, for any formal group
law $F(x,y)$ over a ring~$\Lambda$ there is unique ring homomorphism
$\lambda\:\Omega_U\to\Lambda$ such that
$F(x,y)=\lambda\bigl[\under F(u,v)\bigr]$ (see~\cite{5}). The logarithm of
$\under F(u,v)$ is
$$
\under g(u)=\sum_{n=0}^{\infty}
\frac{\under\varphi[\Bbb CP^n]}{n+1}u^{n+1}=
\sum_{n=0}^{\infty}\frac{\Bbb CP^n}{n+1}u^{n+1}.
$$
Therefore, for any Hirzebruch genus~$\varphi$ we have
$$
g_{\varphi}(u)=\sum_{n=0}^{\infty}\frac{\varphi[\Bbb CP^n]}{n+1}u^{n+1}.
\tag{1}
$$

\noindent 2. The Todd genus $\td(M)$ corresponds to the following series:
$$
Q_{\td}(x)=\frac{x}{1-e^{-x}},\quad
f_{\td}(w)=1-e^{-w},\quad g_{\td}(u)=-\ln (1-u).
$$
These formulae could also be deduced from the fact that
$\td(\Bbb CP^n)=1$. Indeed, using the identity~\thetag{1}, we obtain
$$
g_{\td}(u)=\sum_{n=0}^{\infty}\frac{u^{n+1}}{n+1}=-\ln(1-u).
$$
The corresponding power system is
$$
[u]^{\td}_n=g_{\td}^{-1}\bigl(ng_{\td}(u)\bigr)=1-e^{n\ln(1-u)}=1-(1-u)^n.
$$
Thus, for the Todd genus
$$
f_{\td}(w)=1-e^{-w},\quad g_{\td}(u)=-\ln(1-u),\quad [u]^{\td}_n=1-(1-u)^n.
\tag $2_{\td}$
$$

\noindent 3. The Euler number $e(M)$ (of the tangent bundle). Since
$e(\Bbb CP^n)=n+1$, we deduce (see~\cite{13}) that
$$
g_e(u)=\frac{u}{1-u},\quad f_e(w)=g_e^{-1}(w)=\frac{w}{1+w}.
$$
So $Q(x)=1+x$ and, as might be expected,
$$
e[M]=\biggl(\,\prod_{i=1}^{n}Q(x_i)\biggr)[M]=
(x_1x_2\dots x_n)[M^{2n}]=c_n[M^{2n}].
$$
Next,
$$
[u]_n^e=g_e^{-1}\bigl(ng_e(u)\bigr)=\frac{nu}{1+(n-1)u}.
$$
Thus, for the Euler number,
$$
f_e(w)=\frac{w}{1+w},\quad g_e(u)=\frac{u}{1-u},\quad
[u]^e_n=\frac{nu}{1+(n-1)u}.
\tag $2_e$
$$

\noindent 4. The signature (the~$L$-genus) corresponds
to the following series:
$$
\gather
Q_L(x)=\frac{x}{\th (x)}=\frac{x(e^x+e^{-x})}{e^x-e^{-x}},\quad
f_L(w)=\th(w),
\\
g_L(u)=\arth (u)=\frac12\ln\left(\frac{1+u}{1-u}\right).
\endgather
$$
This is in accordance with Hirzebruch's theorem (the $L$-genus equals
the signature, see~\cite{6}),
from which we deduce that $L(\Bbb CP^{2n})=1$, $L(\Bbb CP^{2n+1})=0$.
The corresponding power system is
$$
[u]^L_n=g_L^{-1}\bigl(ng_L(u)\bigr)=
\frac{e^{\frac n2\ln\frac{1+u}{1-u}}-e^{-\frac n2\ln\frac{1+u}{1-u}}}
{e^{\frac n2\ln\frac{1+u}{1-u}}+e^{-\frac n2\ln\frac{1+u}{1-u}}}=
\frac{(1+u)^n-(1-u)^n}{(1+u)^n+(1-u)^n}.
$$
Thus, for the~$L$-genus,
$$
f_L(w)=\th (w),\quad g_L(u)=\frac12\ln\left(\frac{1+u}{1-u}\right),
\quad [u]^L_n=\frac{(1+u)^n-(1-u)^n}{(1+u)^n+(1-u)^n}.
\tag $2_L$
$$

\noindent 5. There is a one-parametric genus which generalizes three previous
examples,
namely, the $\chi_y$-genus. This is the Hirzebruch genus that
corresponds to the following series:
$$
\gather
Q_{\chi_y}=\frac{x(1+ye^{-x(1+y)})}{1-e^{-x(1+y)}},
\\
f_{\chi_y}(w)=\frac w{Q_{\chi_y}(w)}=
\frac{1-e^{-w(1+y)}}{1+ye^{-w(1+y)}},\quad
g_{\chi_y}(u)=\frac1{1+y}\ln\frac{1+yu}{1-u}.
\endgather
$$
Neglecting the normalizing condition $f_{\chi_y}(w)=w+\cdots$ and replacing
$w(1+y)$ by~$w$, we obtain
$$
\hat f_{\chi_y}(w)=\frac{1-e^{-w}}{1+ye^{-w}},\quad
\hat g_{\chi_y}(u)=\ln\frac{1+yu}{1-u}.
$$
The expression for the power system associated to the
$\chi_y$-genus is the same in both cases:
$$
[u]^{\chi_y}_n=\frac{(1+yu)^n-(1-u)^n}{(1+yu)^n+y(1-u)^n}=
\frac{1-\left(\frac{1-u}{1+yu}\right)^n}
{1+y\left(\frac{1-u}{1+yu}\right)^n}.
\tag $2_{\chi_y}$
$$
For $y=0,-1,1$ we get respectively the Todd genus, the Euler number and
the~$L$-genus.

\smallskip

\noindent 6. The $\hat A$-genus corresponds to the following series:
$$
\gather
Q_A(x)=\frac{x/2}{\sh(x/2)}=\frac{x}{e^{x/2}-e^{-x/2}},
\\
f_A(w)=2\sh\left(\frac{w}2\right),\quad
g_A(u)=2\arsh\left(\frac{u}2\right)
=\ln\left(\frac u2+\sqrt{1+\frac{u^2}{4}}\,\right).
\endgather
$$
The corresponding power system is
$$
\split
[u]^A_n
&=g_A^{-1}\bigl(ng_A(u)\bigr)=2\sh\bigl(n\arsh(u/2)\bigr)\\
&=\left(\frac u2+\sqrt{1+\frac{u^2}{4}}\,\right)^n-
\left(\frac u2+\sqrt{1+\frac{u^2}{4}}\,\right)^{-n}.
\endsplit
$$
Thus, for the $\hat A$-genus,
$$
f_A(w)=2\sh(w/2),\quad g_A(u)=2\arsh(u/2),\quad
[u]^A_n=2\sh\bigl(n\arsh(u/2)\bigr).
\tag $2_A$
$$

\head
\S\,2. Calculation of Hirzebruch genera\\
by means of the Atiyah--Singer index theorem
\endhead

We deal with the following elliptic complexes:
the de Rham complex $\Lambda^*$,
$$
\Gamma(\Lambda^0) @>d>>
\Gamma(\Lambda^1)@>d>>\cdots @>d>> \Gamma(\Lambda^n),
$$
and the Dolbeault complex $\Lambda^{p,*}$,
$$
\Gamma(\Lambda^{p,0})
@>d''>> \Gamma(\Lambda^{p,1}) @>d''>>
\cdots @>d''>> \Gamma(\Lambda^{p,n}),
$$
where $\Lambda^i=\Lambda^i(T^*M)$ is the bundle of differential
$i$-forms on~$M$  nd
$\Lambda^{p,i}=\Lambda^p(T^*M)\wedge\Lambda^i(\ov T^*M)\simeq
\Lambda^p(T^*M)\wedge\Lambda^i(TM)$
is the bundle of differential forms of type $(p,i)$.
(If the manifold~$M$ is endowed with an Hermitian metric, then ${\ov
T}^*M\simeq TM$.)

The following is the Atiyah--Singer index theorem in cohomological form
(see~\cite{3}).

\proclaim{Theorem 2.1}
Let~$M$ be a compact, oriented, differentiable manifold of dimension~$2n$
and let $E=\{d_i\:\Gamma E_i\to\Gamma E_{i+1}\}$ be an
elliptic complex $(i=0,\dots,m-1)$ associated to the tangent bundle
{\rm(}that is, all bundles~$E_i$ are associated to $TM${\rm)}.
Then the index of this complex is determined by the following formula\/{\rm:}
$$
\inda(E)=(-1)^n\Biggl(\biggl(\frac1{e(TM)}\sum_{i=0}^m(-1)^i
\ch(E_i)\biggr)\td(TM\otimes\Bbb C)\Biggr)[M],
$$
where $\ch(E_i)$ is the Chern character of bundle $E_i$.
\endproclaim

\remark{Remark}
Formally factoring the Euler class $e(TM)=x_1\dots x_n$  out of
$\td(TM\otimes\Bbb C)=\prod_{j=1}^n\left(\frac{x_j}{1-e^{-x_j}}
\cdot\frac{-x_j}{1-e^{x_j}}\right)$ in the previous expression, we get the
following formula:
$$
\inda(E)=\Biggl(\biggl(\,\sum_{i=0}^m(-1)^i\ch(E_i)\biggr)
\prod_{j=1}^n\left(\frac{x_j}{1-e^{-x_j}}\cdot
\frac1{1-e^{x_j}}\right)\Biggr)[M].
\tag{3}
$$
\endremark

Let us consider the elliptic complex
$\left\{\Xa_i^y=\bigoplus_{p=0}^ny^p\Lambda^{p,i}\right\}$. (We note that
this object becomes a true elliptic complex only after replacing~$y$ by
actual integers. But its index is obviously defined for arbitrary~$y$ as a
polynomial in~$y$. In what follows we regard such objects as elliptic
complexes.) Using the Atiyah--Singer index theorem, we easily prove the
following fact, which was initially proved by Hirzebruch (see~\cite{6}).

\proclaim{Theorem 2.2}
The index of the elliptic complex $\{\Xa_i^y\}$ equals the
$\chi_y$-genus (defined above) of the manifold~$M$, that is,
$$
\inda(\Xa^y)=\chi_y[M]=\biggl(\,\prod_{j=1}^n
\frac{x_j(1+ye^{-x_j})}{1-e^{-x_j}}\biggr)[M].
$$
\endproclaim

\demo{Proof}
Using formula \thetag{3}, we have:
$$
\split
&\inda(\Xa^y)=\Biggl(\biggl(\,\sum_{i=0}^n(-1)^i\ch(\Xa_i^y)\biggr)
\prod_{j=1}^n\left(\frac{x_j}{1-e^{-x_j}}\frac1{1-e^{x_j}}
\right)\Biggr)[M]
\\
&\qquad=
\Biggl(\ch\biggl(\,\sum_{i=0}^n(-1)^i\Lambda^iTM\biggr)
\ch\biggl(\,\sum_{p=0}^ny^p\Lambda^pT^*M\biggr)\prod_{j=1}^n
\left(\frac{x_j}{1-e^{-x_j}}\frac1{1-e^{x_j}}\right)\Biggr)[M]
\\
&\qquad=
\Biggl(\,\prod_{j=1}^n(1-e^{x_j})\prod_{j=1}^n(1+ye^{-x_j})
\prod_{j=1}^n\left(\frac{x_j}{1-e^{-x_j}}\frac1{1-e^{x_j}}\right)\Biggr)[M]
\\
&\qquad=
\biggl(\,\prod_{j=1}^n\frac{x_j(1+ye^{-x_j})}{1-e^{-x_j}}\biggr)[M]=
\chi_y[M].
\endsplit
$$
Here we have used the following fact (see, for example,~\cite{7}).
\enddemo

\proclaim{Lemma 2.3}
Let $E$ be a complex~$n$-dimensional bundle over a differentiable
mani\-fold~$X$, and let
$c(E)=1+c_1(E)+\cdots+c_n(E)=(1+x_1)\dots \dots(1+x_n)$ be the formal
factorization of the total Chern class. We define
$$
\Lambda_tE:=\sum_{k=0}^{\infty}(\Lambda^kE)t^k,\quad
S_tE:=\sum_{k=0}^{\infty}(S^kE)t^k.
$$
Then
$$
\ch(\Lambda_tE)=\prod_{i=1}^n(1+te^{x_i}),\quad
\ch(S_tE)=\prod_{i=1}^n\frac1{1-te^{x_i}}.
\tag{4}
$$
\endproclaim

{\tolerance=2000
In particular, if we consider the complexes $\{\Lambda^{0,i}\}$,
$\left\{\Cal E_i=\sum_{p=0}^n(-1)^p\Lambda^{p,i}\right\}$,
$\left\{L_i=\sum_{p=0}^n\Lambda^{p,i}\right\}$, we see that their indexes are
respectively $\td[M]=\inda(\Lambda^{0,*})$, $e[M]=\inda(\Cal E)$, and
$L[M]=\inda(L)$.

}
Now we construct an elliptic complex whose index,
under some additional assumptions, equals the
$\hat A$-genus of the manifold~$M$.
Suppose that $c_1(M)\equiv0\mod2$. (This is equivalent to the condition
$w_2(M)=0$, and hence to the existence of a spinor structure on~$M$.) Then
there is a line bundle~$\Cal L$ over~$M$ such that $\Cal L\otimes \Cal
L=\Lambda^nT^*M$, that is, for $c(M)=(1+x_1)\dots(1+x_n)$ we have $c(\Cal
L)=1-\frac{x_1+\cdots+x_n}2$ and $\ch(\Cal
L)=\exp\left(-\frac{x_1+\cdots+x_n}2\right)$. We introduce the complex
$\{A_i=\Lambda^{0,i}\otimes \Cal L\}$.

\proclaim{Theorem 2.4}
The index of the above elliptic complex~$A$ equals the
$\hat A$-genus of~$M$:
$$
\inda(A)=\hat A[M]=
\biggl(\,\prod_{j=1}^n\frac{x_j}{2\sh(x_j/2)}\biggr)[M].
$$
\endproclaim

\demo{Proof}
We again use formulae~\thetag{3} and~\thetag{4}:
$$
\split
&\inda(A)=\Biggl(\ch(\Cal L)\ch\biggl(\,\sum_{i=0}^n(-1)^i\Lambda^{0,i}
\biggr)\prod_{j=1}^n\left(\frac{x_j}{1-e^{-x_j}}
\frac1{1-e^{x_j}}\right)\Biggr)[M]
\\
&\quad=
\Biggl(\exp\left(-\frac{x_1+\cdots+x_n}2\right)
\ch\biggl(\,\sum_{i=0}^n(-1)^i\Lambda^iTM\biggr)
\prod_{j=1}^n\left(\frac{x_j}{1-e^{-x_j}}
\frac1{1-e^{x_j}}\right)\Biggr)[M]
\\
&\quad=
\Biggl(\,\prod_{j=1}^ne^{-x_j/2}\prod_{j=1}^n(1-e^{x_j})
\prod_{j=1}^n\left(\frac{x_j}{1-e^{-x_j}}
\frac1{1-e^{x_j}}\right)\Biggr)[M]
\\
&\quad=
\biggl(\,\prod_{j=1}^n\frac{x_je^{-x_j/2}}{1-e^{-x_j}}\biggr)[M]=
\biggl(\,\prod_{j=1}^n\frac{x_j}{2\sh(x_j/2)}\biggr)[M]=\hat A(M).
\endsplit
$$
\enddemo

Thus we have constructed elliptic complexes associated to the
tangent bundle of~$M$ that enable us to calculate all the
Hirzebruch genera in~\S\,1.

\head
\S\,3. The problem of calculating Hirzebruch genera for manifolds\\
with $\Bbb Z/p$-action in terms of invariants of the action
\endhead

Let $g$ be a transversal endomorphism (that is, $g$ has only finitely many
fixed points) acting on a manifold~$M^{2n}$ such that $g^p=1$ for
some prime $p$. (Hence an action of~$\Bbb Z/p$ is given.) Let $\Cal
P_1,\dots,\Cal P_q$ be the fixed points, and the Jacobi matrix $\Cal J_{\Cal
P_j}(g)$ of the map~$g$ at the point~$\Cal P_j$ ($j=1,\dots,q$)
has eigenvalues
$$
\lambda_k^{(j)}=\exp\biggl(\frac{2\pi i x_k^{(j)}}p\biggr),\quad
x_k^{(j)}\neq 0\mod p,\quad k=1,\dots,n.
$$

Suppose also that there is given an elliptic complex~$E$ on the manifold~$M$.
Let us give the following definition, which is taken from~\cite{1}.

\definition{Definition 3.1}
{\it A lifting of~$g$ to the components of the elliptic complex~$E$} is a
set of linear differential operators
$\varphi_i\:\Gamma(g^*E_i)\to\Gamma(E_i)$.
Here $g^*E_i$ is the pullback of the bundle~$E_i$ under the map~$g$, and
$\Gamma(E_i)$ denotes the linear space of sections of the bundle $E_i$.
\enddefinition

Using~$\varphi_i$, one can define the ``geometric" endomorphisms
$T_i(g,\varphi)\:\Gamma(E_i)\to\Gamma(E_i)$ to be the composite of~$\varphi$
and~$\Gamma_g$, $T_i(g,\varphi)=\varphi_i\circ\Gamma_g$, where
$\Gamma_g\:\Gamma(E_i)\to\Gamma(g^*E_i)$ is the natural map from the sections
of the bundle~$E_i$ into the sections of the pullback $g^*E_i$.  Under these
assumptions we have the following general Atiyah--Bott--Lefschetz theorem
(see~\cite{1}).

{\tolerance=2000
\proclaim{Theorem 3.2}
Suppose that $g\:M\to M$ is a transversal endomorphism of a compact oriented
manifold~$M$. Let $E$ be an elliptic complex on~$M$ and let
$\varphi_i\:\Gamma(g^*E_i)\to\Gamma(E_i)$ be a lifting of~$g$ to the components
of the complex~$E$ such that the corresponding ``geometric" endomorphisms
$T_i(g,\varphi)\:\Gamma(E_i)\to\Gamma(E_i)$ define an endomorphism
$T(g,\varphi)$ of the complex~$E$ {\rm(}that is, the~$T_i$ commute with the
differentials~$d_i${\rm)}. Then the ``equivariant index"
$\inda(g,E):=\sum_{i=1}^n(-1)^{i}\tr T^*_i$ (where
$T_i^*\:H^i(E)\to H^i(E)$) is given by the formula
$$
\inda(g,E)=\sum_{j=1}^q\sigma(\Cal P_j),
\tag{5}
$$
where $\sigma(\Cal P_j)\in\Bbb C$ depends only on local properties of~$T$
and~$\varphi_i$ at the point~$\Cal P_j$.

In particular, if the operator~$\varphi_i$ induces an endomorphism
$\varphi_i(\Cal P_j)\:E_{i,\Cal P_j}\to E_{i,\Cal P_j}$ at each fixed point
$\Cal P_j$, then
$$
\sigma(\Cal P_j)=\sum_{i=0}^n(-1)^i
\frac{\tr\varphi_i(\Cal P_j)}
{\bigl|\det\bigl(1-\Cal J_{\Cal P_j}(g)\bigr)\bigr|}.
\tag{6}
$$
\endproclaim

}
\example{Example 3.1}
If $E_i=\Lambda^i(M)$ is the de Rham complex, then
$\varphi_i(\Cal P_j)=\Lambda^i\Cal J_{\Cal P_j}(g)$.
\endexample

\example{Example 3.2}
If $E_i=\Lambda^{p,i}(M)$ is the Dolbeault complex, then
$$
\varphi_i(\Cal P_j)=\Lambda^p\Cal J'_{\Cal P_j}(g)
\wedge\Lambda^i\barCalJ'_{\Cal P_j}(g),
$$
where $\Cal J'_{\Cal P_j}(g)$ is the holomorphic part of the Jacobi matrix
$\Cal J_{\Cal P_j}(g)$.

Let us find $\sigma(\Cal P_j)$ for the Dolbeault complex $\Lambda^{p,*}$.
Since $\bigl|\det\bigl(1-\Cal J_{\Cal P_j}(g)\bigr)\bigr|=
\det\bigl(1-\Cal J'_{\Cal P_j}(g)\bigr)
\det\bigl(1-\barCalJ'_{\Cal P_j}(g)\bigr)$,
we deduce from \thetag{6} that
$$
\split
\sigma(\Cal P_j)&=\frac{\tr\left(\Lambda^p\Cal J'_{\Cal P_j}(g)
\wedge\sum_i(-1)^i\Lambda^i\barCalJ'_{\Cal P_j}(g)\right)}
{\det\bigl(1-\Cal J'_{\Cal P_j}(g)\bigr)
\det\bigl(1-\barCalJ'_{\Cal P_j}(g)\bigr)}
\\
&=
\frac{\tr\bigl(\Lambda^p\Cal J'_{\Cal P_j}(g)\bigr)
\det\bigl(1-\barCalJ'_{\Cal P_j}(g)\bigr)}
{\det\bigl(1-\Cal J'_{\Cal P_j}(g)\bigr)
\det\bigl(1-\barCalJ'_{\Cal P_j}(g)\bigr)}=
\frac{\tr\Lambda^p\Cal J'_{\Cal P_j}(g)}
{\det\bigl(1-\Cal J'_{\Cal P_j}(g)\bigr)}.
\endsplit
$$
Here we have used the following well-known identity from linear algebra:
$$
\det(1-A)=\sum_i(-1)^i\tr\Lambda^iA
\tag{7}
$$
for any linear operator~$A$. Thus,
$$
\sigma_p(\Cal P_j)=\frac{\tr\Lambda^p\Cal J'_{\Cal P_j}(g)}
{\det\bigl(1-\Cal J'_{\Cal P_j}(g)\bigr)}.
\tag{8}
$$
\endexample

Theorem~3.2 and formula~\thetag{8} enable us to find the fixed point
contribution functions $\sigma(\Cal P_j)$ for the complexes
$\Lambda^{0,*}$, $\Cal E$, $L$,~$\Xa^y$. These complexes calculate
respectively the Todd genus, the Euler number, the~$L$-genus and the
$\chi_y$-genus of the manifold~$M$.

\head
\S\,4. Calculations for the Todd genus,\\
Euler number and the~$L$-genus
\endhead

\subhead 4.1. Calculations for the Euler number \endsubhead
Let us consider the complex $\Cal E_i=\sum_{p=0}^n(-1)^p\Lambda^{p,i}$, for
which $\inda(\Cal E)=e(M)$. In this case
$$
\varphi_i(\Cal P_j)=\sum_{p=0}^n(-1)^p\Lambda^p\Cal J'_{\Cal P_j}(g)
\wedge\Lambda^i\barCalJ'_{\Cal P_j}(g),
$$
and the fixed point contribution functions are as follows
(see~\thetag{8}):
$$
\sigma_e(\Cal P_j)=\sum_{p=0}^n(-1)^p\sigma_p(\Cal P_j)=
\frac{\sum_{p=0}^n(-1)^p\tr\Lambda^p\Cal J'_{\Cal P_j}(g)}
{\det\bigl(1-\Cal J'_{\Cal P_j}(g)\bigr)}=1.
\tag{9}
$$
Here we have used formula \thetag{7}.
From this and theorem~3.2 we deduce that
$$
\inda(g,\Cal E)=\sum_{j=1}^q\sigma_e(\Cal P_j)=q.
$$
Since $\frac1p\sum_{l\in\Bbb Z/p}\inda(g^l,\Cal E)=s$
is the alternating sum of dimensions of the invariant subspaces for the action
of~$g$ on the cohomology of the complex~$\Cal E$, and since
$\inda(1,\Cal E)=\inda(\Cal E)=e(M)$, we have
$$
\inda(1,\Cal E)=e(M)=-\sum_{l=1}^{p-1}\inda(g^l,\Cal E)+ps=
ps-q(p-1)=q+p(s-q).
$$
Therefore, we obtain the following formula for the Euler number:
$$
e(M)\equiv q\pmod p.
\tag{10}
$$

\subhead 4.2. Calculations for the Todd genus \endsubhead
The calculation of the Todd genus of a stably complex
manifold in terms of the
$\Bbb Z/p$-action was carried out by Buchstaber and
Novikov in~\cite{5}. To make our exposition complete, we give
their results here.

\definition{Definition 4.1}
{\it The Atiyah--Bott function} $AB_{\td}(x_1,\dots,x_n)$ of a given fixed
point is the following function of the set of weights
$x_1,\dots,x_n$, \ $x_i\in\Bbb Z/p$:
$$
AB_{\td}(x_1,\dots,x_n)=
-\Tr\Biggl(\,\prod_{k=1}^n\frac1{1-e^{2\pi ix_k/p}}\Biggr),
\tag{11}
$$
where $\Tr\:\Bbb Q(\zeta)\to\Bbb Q$ is the number-theoretical trace,
$\zeta:=e^{2\pi i/p}$.
\enddefinition

It was shown in \cite{5} that
$$
\sum_{j=1}^qAB_{\td}(x_1^{(j)},\dots,x_n^{(j)})\equiv\td(M)\mod p,
\tag{12}
$$
The number-theoretical trace in the definition of Atiyah--Bott functions
for the Todd genus was also calculated in~\cite{5}:
$$
\aligned
AB_{\td}(x_1,\dots,x_n)&\equiv
-\biggl\langle\frac{p[u]^{\td}_{p-1}}{[u]^{\td}_p}
\prod_{k=1}^n\frac u{[u]^{\td}_{x_k}}\biggr\rangle_n\mod p,
\\
AB_{\td}(x_1,\dots,x_n)&\equiv
\sum_{m=0}^n\biggl\langle\frac{pu}{[u]^{\td}_p}\prod_{k=1}^n
\frac u{[u]^{\td}_{x_k}}\biggr\rangle_m\mod p.
\endaligned
\tag{13}
$$
Here $[u]_q^{\td}$ is the $q$th power in the formal group law corresponding
to the Todd genus (see formula~$(2_T)$), and $\bigl\langle h(u)\bigr\rangle_k$
is the coefficient of $u^k$ in the power series $h(u)$.

\subhead 4.3. Calculations for the~$L$-genus \endsubhead
Now let us consider the elliptic complex $L=\sum_{p=0}^n\Lambda^{p,*}$.
Its index is the~$L$-genus: $\inda(L)=L(M)$. In this case, the lifting of~$g$
to the components of the complex~$L$ has the following form:
$$
\varphi_i(\Cal P_j)=\sum_{p=0}^n\Lambda^p\Cal J'_{\Cal P_j}(g)
\wedge\Lambda^i\barCalJ'_{\Cal P_j}(g),
$$
and the fixed point contribution functions are (see formula~\thetag{8})
$$
\split
\sigma_L(\Cal P_j)&=\sum_{p=0}^n\sigma_p(\Cal P_j)=
\frac{\sum_{p=0}^n\tr\Lambda^p\Cal J'_{\Cal P_j}(g)}
{\det\bigl(1-\Cal J'_{\Cal P_j}(g)\bigr)}
\\
&=
\frac{\det\bigl(1+\Cal J'_{\Cal P_j}(g)\bigr)}
{\det\bigl(1-\Cal J'_{\Cal P_j}(g)\bigr)}=
\prod_{k=1}^n\frac{1+e^{2\pi ix_k^{(j)}/p}}{1-e^{2\pi ix_k^{(j)}/p}}.
\endsplit
\tag{14}
$$
Here we have again used formula \thetag{7}.
Hnece it follows from theorem~3.2 that the equivariant index is
$$
\inda(g,L)=\sum_{j=1}^q\sigma_L(\Cal P_j)=
\sum_{j=1}^q\prod_{k=1}^n\frac{1+e^{2\pi ix_k^{(j)}/p}}
{1-e^{2\pi ix_k^{(j)}/p}}.
$$
As before, $\frac1p\sum_{l\in\Bbb Z/p}\inda(g^l,L)=s$ is the alternating sum
of the dimensions of the invariant subspaces
for the action of~$g$ on the cohomology of the complex~$L$, and
we have $\inda(1,L)=\inda(L)=L(M)$. Hence
$$
L(M)=\inda(1,L)=-\sum_{l=1}^{p-1}\inda(g^l,L)+ps=
-\sum_{j=1}^q\sum_{l=1}^{p-1}\prod_{k=1}^n
\frac{1+e^{2\pi ix_k^{(j)}l/p}}{1-e^{2\pi ix_k^{(j)}l/p}}+ps.
$$
Again, we consider the number-theoretical trace
$\Tr\:\Bbb Q(\zeta)\to\Bbb Q$
and introduce the Atiyah--Bott functions $AB_L(x_1,\dots,x_n)$ as
$$
AB_L(x_1,\dots,x_n)=-\Tr\Biggl(\,\prod_{k=1}^n
\frac{1+e^{2\pi ix_k/p}}{1-e^{2\pi ix_k/p}}\Biggr).
\tag{15}
$$
Then we have
$$
\sum_{j=1}^qAB_L\bigl(x_1^{(j)},\dots,x_n^{(j)}\bigr)\equiv L(M)\mod p.
\tag{16}
$$
Relations \thetag{15} and \thetag{16} are analogous to relations
\thetag{11},~\thetag{12} for the Todd genus. It remains
to calculate the number-theoretical trace in the definition
of the Atiyah--Bott functions $AB_L(x_1,\dots,x_n)$.

We set ${\theta=\frac{1-\zeta}{1+\zeta}}$.
Then ${\zeta=e^{2\pi i/p}=\frac{1-\theta}{1+\theta}}$.
We must calculate
$\Tr\left(\prod_{k=1}^n\frac{1+\zeta^{x_k}}{1-\zeta^{x_k}}\right)$.
We shall carry out our calculations in the~$p$-adic extension
$\Bbb Q_p(\zeta)$
of the field $\Bbb Q(\zeta)$. Let us write $A\simeq B$ if
$A$ and $B$ are equal modulo $p\Bbb Z_p\subset\Bbb Q_p$. First of all, we
prove the following statement.

\proclaim{Lemma 4.2}
$\Tr(\theta^k)\simeq0$ for any $k\geq1$.
\endproclaim

\demo{Proof}
Since $2=\zeta^p+1=(\zeta+1)(\zeta^{p-1}-\zeta^{p-2}+
\zeta^{p-3}-\cdots-\zeta+1)$, we have
$$
\theta=\frac{(1-\zeta)
(\zeta^{p-1}-\zeta^{p-2}+\zeta^{p-3}-\cdots-\zeta+1)}2=
\zeta^{p-1}-\zeta^{p-2}+\zeta^{p-3}-\cdots-\zeta.
$$
Therefore,
$$
\Tr\theta^k=\sum_{m=1}^{p-1}\bigl((\zeta^m)^{p-1}-(\zeta^m)^{p-2}+
(\zeta^m)^{p-3}-\cdots-\zeta^m\bigr)^k\simeq0
$$
because
$$
\sum_{m=1}^{p-1}\bigl((\zeta^m)^{p-1}-(\zeta^m)^{p-2}+\cdots-
\zeta^m\bigr)^k=\sum_{m=0}^{p-1}\bigl((\zeta^m)^{p-1}-
(\zeta^m)^{p-2}+\cdots-\zeta^m\bigr)^k,
$$
and because $\sum_{m=0}^{p-1}(\zeta^m)^r\simeq0$ for any~$r$.
The lemma is proved.
\enddemo

We further deduce that
$$
\prod_{k=1}^n\frac{1+\zeta^{x_k}}{1-\zeta^{x_k}}=
\prod_{k=1}^n\frac{1+\left(\frac{1-\theta}{1+\theta}\right)^{x_k}}
{1-\left(\frac{1-\theta}{1+\theta}\right)^{x_k}}=
\frac1{\theta^n}\prod_{k=1}^n
\frac{\theta\left(1+\left(\frac{1-\theta}{1+\theta}\right)^{x_k}\right)}
{1-\left(\frac{1-\theta}{1+\theta}\right)^{x_k}}=:
\frac1{\theta^n}\sum_{k=0}^{\infty}A_k\theta^k,
$$
where the $A_k\in\Bbb Z_p$ are $p$-adic integers
(since $1/x_k\in\Bbb Z_p$). Therefore, we can write
$$
\Tr\biggl(\frac1{\theta^n}\sum_{k=0}^{\infty}A_k\theta^k\biggr)\simeq
\Tr\biggl(\frac1{\theta^n}\sum_{k=0}^nA_k\theta^k\biggr)=
\sum_{k=0}^n\bigl(A_k\Tr(\theta^{k-n})\bigr).
$$
Set $\Tr\theta^{-s}=B_s$ and introduce two formal power series
$A(u)=\sum_{k=0}^{\infty}A_ku^k$, \ $B(u)=\sum_{k=0}^{\infty}B_ku^k$.
It follows from~$(2_L)$ that
$A(u)=\prod_{k=1}^n\frac{u}{[u]_{x_k}^L}$.
Therefore, we must calculate the coefficient of~$u^n$ in the series
$A(u)B(u)$. We have
$$
\split
B(u)&=\Tr\biggl(1+\sum_{s=1}^{\infty}\theta^{-s}u^s\biggr)=
\Tr\left(\frac1{1-\theta^{-1}u}\right)
\\
&=
\Tr\left(\frac{\theta}{\theta-u}\right)=
\Tr\left(1+\frac u{\theta-u}\right)=
(p-1)+u\Tr\left(\frac1{\theta-u}\right).
\endsplit
$$

Observe that if $\varphi_{\alpha}(u)$ is the minimal polynomial for an
element~$\alpha$ with respect to the extension $\Bbb Q_p(\zeta)\mid\Bbb Q_p$,
then $\Tr\frac1{\alpha-u}=
-\frac{\varphi_{\alpha}'(u)}{\varphi_{\alpha}(u)}$. Since
$$
\split
0&=\frac{\zeta^p-1}{\zeta-1}=\frac{\left(
\frac{1-\theta}{1+\theta}\right)^p-1}
{\left(\frac{1-\theta}{1+\theta}\right)-1}
\\
&=
\frac{\bigl((1-\theta)^p-(1+\theta)^p\bigr)(1+\theta)}
{(1+\theta)^p(-2\theta)}=\frac1{(1+\theta)^{p-1}}
\frac{(1+\theta)^p-(1-\theta)^p}{2\theta},
\endsplit
$$
we deduce that $\varphi_{\theta}(u)=\frac{(1+u)^p-(1-u)^p}{2u}$ is the
minimal polynomial for $\theta=\frac{1-\zeta}{1+\zeta}$. Hence,
$$
\gather
-\Tr\frac1{\theta-u}=\frac{\varphi_{\theta}'(u)}{\varphi_{\theta}(u)}=
p\frac{(1+u)^{p-1}+(1-u)^{p-1}}{(1+u)^p-(1-u)^p}-\frac1u,
\\
B(u)=(p-1)+u\Tr\left(\frac1{\theta-u}\right)=
\frac{(1+u)^{p-1}-(1-u)^{p-1}}{(1+u)^p-(1-u)^p},
\\
\Tr\biggl(\,\prod_{k=1}^n\frac{1+\zeta^{x_k}}{1-\zeta^{x_k}}\biggr)\simeq
\bigl\langle A(u)B(u)\bigr\rangle_n=
\biggl\langle p\frac{(1+u)^{p-1}-(1-u)^{p-1}}{(1+u)^p-(1-u)^p}
\prod_{k=1}^n\frac u{[u]_{x_k}^L}\biggr\rangle_n,
\endgather
$$
We deduce that
$$
AB_L(x_1,\dots,x_n)\equiv-\biggl\langle
p\frac{(1+u)^{p-1}-(1-u)^{p-1}}{(1+u)^p-(1-u)^p}
\prod_{k=1}^n\frac u{[u]_{x_k}^L}\biggr\rangle_n\mod p.
$$
This formula is analogous to the first formula in \thetag{13}. Further,
$$
\split
&p\frac{(1+u)^{p-1}-(1-u)^{p-1}}{(1+u)^p-(1-u)^p}=
p\frac{(1-u)(1+u)^p-(1+u)(1-u)^p}{(1-u)(1+u)\bigl((1+u)^p-(1-u)^p\bigr)}
\\
&\qquad=
p\frac{(1+u)^p-(1-u)^p-u\bigl((1+u)^p+(1-u)^p\bigr)}
{(1-u^2)\bigl((1+u)^p-(1-u)^p\bigr)}
\\
&\qquad=
\frac p{1-u^2}-\frac{pu}{(1-u^2)
\frac{(1+u)^p-(1-u)^p}{(1+u)^p+(1-u)^p}}=
\frac p{1-u^2}-\frac{pu}{(1-u^2)[u]_p^L}
\\
&\qquad\simeq
-\frac{pu}{[u]_p^L}(1+u^2+u^4+\cdots).
\endsplit
$$
From this we obtain
$$
\split
AB_L(x_1,\dots,x_n)&\equiv
-\biggl\langle p\frac{(1+u)^{p-1}-(1-u)^{p-1}}{(1+u)^p-(1-u)^p}
\prod_{k=1}^n\frac u{[u]_{x_k}^L}\biggr\rangle_n
\\
&\equiv
\biggl\langle \frac{pu}{[u]_p^L}\prod_{k=1}^n\frac u{[u]^L_{x_k}}
(1+u^2+u^4+\cdots)\biggr\rangle_n\mod p,
\endsplit
$$
whence
$$
AB_L(x_1,\dots,x_n)\equiv\sum_{i=0}^{[n/2]}\biggl\langle
\frac{pu}{[u]_p^L}\prod_{k=1}^n\frac u{[u]^L_{x_k}}\biggr\rangle_{n-2i}.
$$
{\tolerance=2000
At the same time,
$[u]^L_k=\frac{(1+u)^k-(1-u)^k}{(1+u)^k+(1-u)^k}=\th(k\arth u)$,
and the series~$\arth u$
(as well as $\th u$) contains only odd powers
of~$u$. Therefore, the series
$\frac{pu}{[u]_p^L}\prod_{k=1}^n\frac u{[u]_{x_k}^L}$
contains only even powers of~$u$. Thus we can finally write the formulae
for the Atiyah--Bott functions for the~$L$-genus which are analogous to
formulae~\thetag{13} for the Todd genus:
$$
\aligned
AB_L(x_1,\dots,x_n)&\equiv-\biggl\langle
p\frac{(1+u)^{p-1}-(1-u)^{p-1}}{(1+u)^p-(1-u)^p}
\prod_{k=1}^n\frac u{[u]_{x_k}^L}\biggr\rangle_n\mod p,
\\
AB_L(x_1,\dots,x_n)&\equiv\sum_{m=0}^n\biggl\langle
\frac{pu}{[u]^L_p}\prod_{k=1}^n\frac u{[u]_{x_k}^L}\biggr\rangle_m\mod p.
\endaligned
\tag{17}
$$
By \thetag{16}, we see that
$$
L(M)\equiv\sum_{j=1}^q\sum_{m=0}^n\biggl\langle\frac{pu}{[u]_p^L}
\prod_{k=1}^n\frac u{[u]_{x_k^{(j)}}^L}\biggr\rangle_m\mod p,
\tag 18
$$
where $[u]_p^L=\frac{(1+u)^p-(1-u)^p}{(1+u)^p+(1-u)^p}$.

}
\remark{Remark}
Relations similar to \thetag{12}, \thetag{16} could also be
obtained for the Euler number. In this case, the Atiyah--Bott functions are
$$
AB_e(x_1,\dots,x_n)=
-\Tr\biggl(\,\prod_{k=1}^n\frac{1-\zeta}{1-\zeta}\biggr)=
-\Tr(1)=-(p-1)\equiv1\mod p.
$$
Therefore,
$$
q=\sum_{j=1}^qAB_e(x_1^{(j)},\dots,x_n^{(j)})\equiv e(M)\mod p.
\tag{19}
$$
This relation is analogous to relations \thetag{12} and \thetag{16} for
the Todd genus and the~$L$-genus.
\endremark

\head
\S\,5. General results on the calculation of\\
Hirzebruch genera via invariants of the $\Bbb Z/p$-action \endhead

{\tolerance=2000
Here we consider another approach to calculating the equivariant index
$\inda(g,E)=\sum_{i=0}^m(-1)^i\tr(g,H^i)$ of an elliptic complex~$E$.
This approach is taken from~\cite{3} (see also~\cite{7}).

}
In what follows we adopt somewhat weaker assumptions about the action
of an operator~$g$, $g^p=1$, on a stably complex manifold~$M^{2n}$. Namely,
we remove the transversality condition. Let
$M^g=\bigl\{x\in M\mid gx=x\bigr\}$ be the fixed point set, and let
$M^g=\cup M_{\nu}^g$ be its decomposition into connected components. Then
the equivariant index can be computed as the sum of the contributions
$\sigma(M_{\nu}^g)$ corresponding to the fixed point components $M_{\nu}^g$
(see~\cite{7}). These contributions are calculated as follows.

Let $Y=M_{\nu}^g$ be one of the fixed point components of $M^{2n}$.
For each point $p\in Y$, $g$ acts linearly on the tangent space
$T_pM$. This tangent space decomposes into the direct sum of the
eigenspaces~$N_{p,\lambda}$ for eigenvalues~$\lambda$, $|\lambda|=1$.
In this way we obtain the eigenbundle~$N_{\lambda}$ over~$Y$. The
$N_1$ is just the tangent bundle $TY$ to~$Y$. With
$d_{\lambda}=\rk N_{\lambda}$, we therefore have:
$$
TM|_Y=\bigoplus_{\lambda}N_{\lambda},\qquad
c(N_{\lambda})=\prod_{i=1}^{d_{\lambda}}(1+x_i^{\lambda}),
\tag{20}
$$
that is,
$$
c(TM_Y)=\prod_{\lambda}\prod_{i=1}^{d_{\lambda}}(1+x_i^{\lambda})=
\prod_{i=1}^n(1+x_i).
$$
The recipe for calculating $\sigma(Y)$ is as follows.
Consider the index formula~\thetag{3} in Theorem~2.1:
$$
\inda(E)=\Biggl(\biggl(\,\sum_{i=0}^m(-1)^i\ch(E_i)\biggr)c_n(M)
\prod_{j=1}^n\left(\frac1{1-e^{-x_j}} \frac1{1-e^{x_j}}\right)\Biggr)[M]
\tag{21}
$$
and replace~$M$ by~$Y$ and $e^{x_i}$ by $\lambda^{-1}e^{x_i}$
(where $x_i$ belongs to the eigenvalue $\lambda$). Apply the same process to
the terms~$\ch(E_i)$. This can obviously be done if
the~$E_i$ are associated to the tangent bundle of~$M$. This
``recipe" is taken from~\cite{7}.

In the case of finite number of fixed points we have $Y=\bold p\bold t$, \
$c_n(Y)=1$, and so we must replace $x_1x_2\dots x_n$ by~1, and
$e^{x_j^{\lambda}}$ by~$\lambda_j^{-1}$. Therefore, introducing
the ``weight"~$x_i$ by the
formula~$\lambda_j=\exp\left(\frac{2\pi ix_j}p\right)$,
we just have to replace~$x_j^{\lambda}$ by $-\frac{2\pi i}px_j$. (Here
in the first case~$x_j^{\lambda}$ stands for the first Chern class of
``virtual" line subbundle in $TM$ corresponding to the eigenvalue~$\lambda_j$,
while in the second case $x_j$ is the ``weight" of the fixed point and
is defined only modulo~$p$.)

\example{Example 5.1}
Let us consider the $\chi_y$-genus of a manifold~$M$. Applying our
recipe to the formula from Theorem~2.2, we obtain the following formula for
the contribution of each fixed point~$\Cal P$:
$$
\sigma(\Cal P)=\prod_{k=1}^n\frac{1+ye^{2\pi ix_k/p}}
{1-e^{2\pi ix_k/p}}.
\tag{22}
$$
Putting $y=-1,0,1$, we obtain the formulae for the contribution function for
the Euler number, the Todd genus and the~$L$-genus:
$$
\sigma_e(\Cal P)=1,\qquad
\sigma_{\td}(\Cal P)=\prod_{k=1}^n\frac1{1-e^{2\pi ix_k/p}},\qquad
\sigma_L(\Cal P)=\prod_{k=1}^n\frac{1+e^{2\pi ix_k/p}}
{1-e^{2\pi ix_k/p}}.
$$
These formulae coincide with \thetag{9}, \thetag{14}
and the formula in~\cite{5}, which were deduced from the Atiyah--Bott
theorem~3.2.
\endexample

Now consider the general case of an arbitrary Hirzebruch
genus~$\varphi$:
$$
\varphi(M)=\biggl(\,\prod_{i=1}^n\frac{x_i}{f_{\varphi}(x_i)}\biggr)[M],
$$
where $g_{\varphi}(u)=f_{\varphi}^{-1}(u)$ is the logarithm of the
corresponding formal group law. Suppose that there is an elliptic complex
$E_{\varphi}$ associated to $TM$ whose index equals $\varphi(M)$.
Applying the above recipe, we see that the contribution functions of the
fixed points for the action of~$g$ on~$M$ are given by the formula
$$
\sigma(\Cal P)=\prod_{k=1}^n\frac1{f_{\varphi}(-2\pi ix_k/p)}.
\tag{23}
$$
Here the $x_k$ are the ``weights" of the fixed point $\Cal P$.
They are determined by the formula
$\lambda_k=\exp\left(\frac{2\pi ix_k}p\right)$, \ $x_k\neq0\mod p$, where the
$\lambda_k$ are the eigenvalues of the Jacobi matrix $\Cal J_{\Cal P}(g)$
of the map~$g$ at the point~$\Cal P$. The equivariant index of~$E_{\varphi}$
is then determined by the formula
$$
\inda(g,E_{\varphi})=\sum_{j=1}^q\prod_{k=1}^n
\frac 1{f_{\varphi}\bigl(-2\pi ix_k^{(j)}/p\bigr)}.
$$

Further, $\frac1p\sum_{l\in\Bbb Z/p}\inda(g^l,E_{\varphi})=s$ is the
alternating sum of the dimensions of the equivariant subspaces for the action
of~$g$ on the cohomology of~$E_{\varphi}$, and
$\inda(1,E_{\varphi})=\inda(E_{\varphi})=\varphi(M)$. Therefore,
$$
\varphi(M)=\inda(1,E_{\varphi})=-\sum_{l=1}^{p-1}
\inda(g^l,E_{\varphi})+ps=-\sum_{j=1}^q\sum_{l=1}^{p-1}\prod_{k=1}^n
\frac1{f_{\varphi}\bigl(-2\pi ix_k^{(j)}l/p\bigr)}+ps.
$$
Consider the number-theoretical trace
$\Tr\:\Bbb Q_p(\zeta)\to\Bbb Q_p$, where
$\zeta:=\exp\frac{2\pi i}p$. Then
$$
\sum_{l=1}^{p-1}\prod_{k=1}^n
\frac1{f_{\varphi}\bigl(-2\pi ix_k^{(j)}l/p\bigr)}=
\Tr\prod_{k=1}^n\frac1{f_{\varphi}\bigl(-2\pi ix_k^{(j)}/p\bigr)}.
$$

Thus, we obtain the following result.

\proclaim{Theorem 5.1}
Suppose that there is an elliptic complex of bundles associated to $TM$
whose index is equal to the Hirzebruch genus $\varphi(M)$ of the manifold
$M$. Let~$g$ be a holomorphic transversal endomorphism
acting on~$M$ such that $g^p=1$. Then we have the following formula
for~$\varphi(M)${\rm:}
$$
\varphi(M)\equiv-\sum_{j=1}^q\Tr\prod_{k=1}^n
\frac1{f_{\varphi}\bigl(-2\pi ix_k^{(j)}/p\bigr)}\mod p.
$$
\endproclaim

It is now convenient to make the following definition.

\definition{Definition 5.2}
{\it The Atiyah--Bott fixed point function $AB_{\varphi}(x_1,\dots,x_n)$
corresponding to the genus}~$\varphi$ is defined to be the following function
of the set of weights $x_1,\dots,x_n$, \ $x_i\in\Bbb Z/p$:
$$
AB_{\varphi}(x_1,\dots,x_n)=-\Tr\prod_{k=1}^n
\frac1{f_{\varphi}\bigl(-2\pi ix_k/p\bigr)}.
$$
\enddefinition

{\tolerance=2000
Then we immediately get
$$
\sum_{j=1}^qAB_{\varphi}\bigl(x_1^{(j)},\dots,x_n^{(j)}\bigr)\equiv
\varphi(M)\mod p.
\tag{24}
$$
Now we set $\theta=f_{\varphi}(-2\pi i/p)$. Then
$-2\pi i/p=g_{\varphi}(\theta)$, and $f_{\varphi}(-2\pi i/px_k)=
f_{\varphi}\bigl(x_kg_{\varphi}(\theta)\bigr)=
[\theta]_{x_k}^{\varphi}$. Hence the following statement holds:

}
\proclaim{Proposition 5.3}
The Atiyah--Bott fixed point function $AB_{\varphi}(x_1,\dots,x_n)$
corresponding to the genus~$\varphi$ can be computed as
$$
AB_{\varphi}(x_1,\dots,x_n)=-\Tr\biggl(\,\prod_{k=1}^n
\frac1{[\theta]_{x_k}^{\varphi}}\biggr),\qquad
\theta=f_{\varphi}\left(\frac{-2\pi i}p\right).
$$
\endproclaim

\proclaim{Lemma 5.4}
For any $k>0$ we have $\Tr\theta^k\simeq0$
{\rm(}that is, $\Tr\theta^k\in p\Bbb Z_p${\rm)}.
\endproclaim

\demo{Proof}
Let $\theta=f_{\varphi}\left(-\frac{2\pi i}p\right)=
C_0+C_1\zeta+\cdots+C_{p-1}\zeta^{p-1}\in\Bbb Q_p(\zeta)$, \
$C_i\in\Bbb Z_p$. Then $C_0+C_1\zeta^m+C_2\zeta^{2m}+\dots+
C_{p-1}\zeta^{(p-1)m}=f_{\varphi}\left(-\frac{2\pi i}pm\right)$.
In particular, $C_0+C_1+\cdots+C_{p-1}=0$.
Therefore,
$$
\split
\Tr\theta^k&=\sum_{m=1}^{p-1}f_{\varphi}\left(-\frac{2\pi i}pm\right)=
\sum_{m=1}^{p-1}\bigl(C_0+C_1\zeta^m+C_2\zeta^{2m}+\dots+
C_{p-1}\zeta^{(p-1)m}\bigr)
\\
&=
\sum_{m=0}^{p-1}\bigl(C_0+C_1\zeta^m+C_2\zeta^{2m}+\dots+
C_{p-1}\zeta^{(p-1)m}\bigr)\simeq0,
\endsplit
$$
since $\sum_{m=0}^{p-1}\zeta^{rm}\simeq0$ for any~$r$.
The lemma is proved.
\enddemo

\example{Example 5.2}
For the Todd genus, $\theta=f_{\td}\left(-\frac{2\pi i}p\right)=
1-e^{2\pi i/p}=1-\zeta$ (see formula~$(2_{\td})$). Hence,
$C_0=1$, \ $C_1=-1,C_i=0$ for $i>1$.
\endexample

\example{Example 5.3}
For the~$L$-genus, $\theta=\hat f_L\left(-\frac{2\pi i}p\right)=
\frac{1-e^{2\pi i/p}}{1+e^{2\pi i/p}}=\frac{1-\zeta}{1+\zeta}=
\zeta^{p-1}-\zeta^{p-2}+\zeta^{p-3}-\cdots-\zeta$
(see Lemma~4.2). Hence, $C_0=0$, \ $C_{2i+1}=-1$, \ $C_{2i}=1$ for $i>1$.
\endexample

\head
\S\,6. Calculations for the $\hat A$-genus and the $\chi_y$-genus
\endhead

\subhead 6.1. Calculations for the $\hat A$-genus \endsubhead
We consider the $\hat A$-genus
$$
\hat A(M)=\biggl(\,\prod_{j=1}^n
\frac{x_je^{-x_j/2}}{1-e^{-x_j}}\biggr)[M]
$$
for a stably complex manifold~$M$ such that
$c_1(M)\equiv0\mod2$. In Theorem~2.4, we gave an elliptic complex whose index
is $\hat A(M)$. The Atiyah--Bott functions for the
$\hat A$-genus are as follows (see Definition~5.2):
$$
\gather
AB_{\hat A}(x_1,\dots,x_n)=-\Tr\left(\,\prod_{k=1}^n
\frac{\exp\left(\frac{2\pi ix_k}{2p}\right)}
{1-\exp\left(\frac{2\pi ix_k}p\right)}\right)=-\Tr\biggl(\,\prod_{k=1}^n
\frac{\zeta^{x_k/2}}{1-\zeta^{x_k}}\biggr),
\tag{25}
\\
\sum_{j=1}^qAB_{\hat A}\bigl(x_1^{(j)},\dots,x_n^{(j)}\bigr)\equiv
\hat A(M)\mod p.
\tag{26}
\endgather
$$

Now we have to calculate the number-theoretical trace in the definition of
the Atiyah--Bott functions $AB_{\hat A}(x_1,\dots,x_n)$. By
Proposition~5.3,
$\frac{\zeta^{x_k/2}}{1-\zeta^{x_k}}=\frac1{[\theta]_{x_k}^A}$, where
$[u]_m^A=2\sh\left(m\arsh\left(\frac u2\right)\right)$ is the $m$th power in
the formal group law corresponding to the $\hat A$-genus (see
formula~$(2_A)$), and $\theta=2\sh\left(-\frac{2\pi i}{2p}\right)=
-2\sh\left(\frac{\pi i}p\right)=e^{-\pi i/p}-e^{\pi i/p}=
\zeta^{-1/2}-\zeta^{1/2}=\zeta^{(p+1)/2}-\zeta^{(p-1)/2}$.

We claim that the minimal polynomial for the element
$\theta=-2\sh\left(\frac{\pi i}p\right)$ is
$$
\varphi_{\theta}(u)=\frac{2\sh\left(p\arsh\left(\frac u2\right)\right)}u.
$$
Indeed,
$$
\varphi_{\theta}(\theta)=
\frac{2\sh\left(p\arsh\left(\left(-2\sh\left(
\frac{\pi i}p\right)\right)\!\Big/2\right)\right)}{\theta}=
\frac{-2\sh\left(p\cdot\frac{\pi i}p\right)}\theta=0
$$
and $\varphi_{\theta}(u)$ is a polynomial of degree $p-1$ with
leading term~$u^{p-1}$. For instance, for $p=3$ we get
$\varphi_{\theta}(u)=u^2+3$, for $p=5$ we get
$\varphi_{\theta}(u)=u^4+5u^2+5$, and so on.

Further, we have
$$
\split
\prod_{k=1}^n\frac{\zeta^{x_k/2}}{1-\zeta^{x_k}}&=
\prod_{k=1}^n\frac1{f_{\hat A}\bigl(x_kg_{\hat A}(\theta)\bigr)}=
\prod_{k=1}^n\frac1{[\theta]^{\hat A}_{x_k}}
\\
&=
\frac1{\theta^n}\prod_{k=1}^n\frac{\theta}{[\theta]^{\hat A}_{x_k}}=:
\frac1{\theta^n}\sum_{i=0}^{\infty}A_i\theta^i=\frac1{\theta^n}A(\theta),
\endsplit
$$
\pagebreak

\noindent where $A(u):=\prod_{k=1}^n\frac{u}{[u]^{\hat A}_{x_k}}=
\sum_{i=0}^{\infty}A_iu^i$, \ $A_i\in\Bbb Z_p$ are $p$-adic integers.
It follows from lemma~5.4 that $\Tr(\theta^k)\simeq0$ for any $k>0$.
(Here $\theta=\zeta^{(p+1)/2}-\zeta^{(p-1)/2}$.) So the required
trace is
$$
\split
\Tr\biggl(\,\prod_{k=1}^n\frac{\zeta^{x_k/2}}{1-\zeta^{x_k}}\biggr)&=
\Tr\biggl(\frac1{\theta^n}\sum_{i=0}^{\infty}A_i\theta^i\biggr)\simeq
\Tr\biggl(\frac1{\theta^n}\sum_{i=0}^nA_i\theta^i\biggr)
\\
&=
\sum_{i=0}^n\bigl(A_i\Tr(\theta^{i-n})\bigr)=
\bigl\langle A(u)B(u)\bigr\rangle_n,
\endsplit
$$
where $B_s:=\Tr\theta^{-s}$, \ $B(u):=\sum_{i=0}^{\infty}B_iu^i$.
As before, $B(u)=(p-1)+u\Tr\left(\frac1{\theta-u}\right)$, \
$\Tr\left(\frac1{\theta-u}\right)=
-\frac{\varphi'_{\theta}(u)}{\varphi_{\theta}(u)}$.
Therefore,
$$
\split
\varphi_{\theta}(u) &=
\frac{2\sh\left(p\arsh\left(\frac u2\right)\right)}u,
\\
\varphi'_{\theta}(u) &=
\frac{p\cosh\left(p\arsh\left(\frac u2\right)\right)}
{u\sqrt{1+\frac{u^2}4}}-
\frac{2\sh\left(p\arsh\left(\frac u2\right)\right)}{u^2},
\\
\frac{\varphi'_{\theta}(u)}{\varphi_{\theta}(u)} &=
\frac{p}{2\sqrt{1+\frac{u^2}4}}
\frac{\cosh\left(p\arsh\left(\frac u2\right)\right)}
{\sh\left(p\arsh\left(\frac u2\right)\right)}-\frac1u,
\endsplit
$$
$$
\multline
B(u)=p-1-u\frac{\varphi'_{\theta}(u)}{\varphi_{\theta}(u)}=
p-\frac{p\frac u2}{\sqrt{1+\frac{u^2}4}}
\frac{\cosh\left(p\arsh\left(\frac u2\right)\right)}
{\sh\left(p\arsh\left(\frac u2\right)\right)}
\\
=
\frac{p\sh\left((p-1)\arsh\left(\frac u2\right)\right)}
{\cosh\left(\arsh\left(\frac u2\right)\right)
\sh\left(p\arsh\left(\frac u2\right)\right)}=
\frac{p\sh\left((p-1)\arsh\left(\frac u2\right)\right)}
{\sqrt{1+\frac{u^2}4}\,\sh\left(p\arsh\left(\frac u2\right)\right)}.
\endmultline
$$
Here we have used the following formulae: $\cosh\bigl(\arsh(u/2)\bigr)=
\sqrt{1+u^2/4}$ and $\sh(x-y)=\sh x\cosh y-\sh y\cosh x$.
We further deduce that
$$
\split
B(u)&=\frac{p\sh\left((p-1)\arsh\left(\frac u2\right)\right)}
{\cosh\left(\arsh\left(\frac u2\right)\right)
\sh\left(p\arsh\left(\frac u2\right)\right)}
\\
&=
\frac{2p\sh\left((p-1)\arsh\left(\frac u2\right)\right)}
{\sh\left((p-1)\arsh\left(\frac u2\right)\right)+
\sh\left((p+1)\arsh\left(\frac u2\right)\right)}
\\
&=
p\frac{2[u]^A_{p-1}}{[u]^A_{p-1}+[u]^A_{p+1}}\simeq
p\frac{[u]^A_{p-1}-[u]^A_{p+1}}{[u]^A_{p-1}+[u]^A_{p+1}}.
\endsplit
$$
Here $[u]^A_m=2\sh\bigl(m\arsh(u/2)\bigr)$, and we have used the formula
$2\cosh x\sh y=\sh(y+x)+\sh(y-x)$. Thus we obtain the following formulae
for the Atiyah--Bott functions $AB_{\hat A}(x_1,\dots,x_n)= -\Tr\left(\,
\prod_{k=1}^n\frac{\zeta^{x_k/2}}{1-\zeta^{x_k}}\right)$
for the $\hat A$-genus:
$$
\aligned
AB_{\hat A}(x_1,\dots,x_n)&=-\biggl\langle
p\frac{\sh\left((p-1)\arsh\left(\frac u2\right)\right)}
{\sqrt{1+\frac{u^2}4}\,\sh\left(p\arsh\left(\frac u2\right)\right)}
\prod_{k=1}^n\frac{u}{[u]^{\hat A}_{x_k}}\biggr\rangle_n,
\\
AB_{\hat A}(x_1,\dots,x_n)&=-\biggl\langle
p\frac{2[u]^A_{p-1}}{[u]^A_{p-1}+[u]^A_{p+1}}
\prod_{k=1}^n\frac{u}{[u]^{\hat A}_{x_k}}\biggr\rangle_n,
\\
AB_{\hat A}(x_1,\dots,x_n)&\simeq\biggl\langle
p\frac{[u]^A_{p+1}-[u]^A_{p-1}}{[u]^A_{p-1}+[u]^A_{p+1}}
\prod_{k=1}^n\frac{u}{[u]^{\hat A}_{x_k}}\biggr\rangle_n.
\endaligned
\tag{27}
$$
Then the $\hat A$-genus itself is calculated as $\hat A(M)\equiv
\sum_{j=1}^qAB_{\hat A}\bigl(x_1^{(j)},\dots,x_n^{(j)}\bigr)$.

\subhead 6.2. Calculations for the $\chi_y$-genus \endsubhead
Let us consider the $\chi_y$-genus
$$
\chi_y(M)=\biggl(\,\prod_{j=1}^n
\frac{x_j(1+ye^{-x_j})}{1-e^{-x_j}}\biggr)[M]
$$
of a manifold~$M$. An elliptic complex whose index is
$\chi_y(M)$ is given by Theorem~2.2. The Atiyah--Bott functions for
the $\chi_y(M)$-genus are as follows (see Definition~5.2):
$$
\gather
{\align
AB_{\chi_y}(x_1,\dots,x_n)&=-\Tr\left(\,\prod_{k=1}^n
\frac{1+y\exp\left(\frac{2\pi ix_k}{p}\right)}
{1-\exp\left(\frac{2\pi ix_k}p\right)}\right)
\\
&=-\Tr\biggl(\,
\prod_{k=1}^n\frac{1+y\zeta^{x_k}}{1-\zeta^{x_k}}\biggr),
\tag{28}
\endalign}
\\
\sum_{j=1}^qAB_{\chi_y}\bigl(x_1^{(j)},\dots,x_n^{(j)}\bigr)\equiv
\chi_y(M)\mod p.
\tag{29}
\endgather
$$
These formulae could be also deduced from the Atiyah--Bott theorem~3.2
as was done for the $L$-genus in~\S\,4.3.

We shall now calculate the number-theoretical trace in the definition
of the Atiyah--Bott functions $AB_{\chi_y}(x_1,\dots,x_n)$. By
Proposition~5.3,
$\frac{1+y\zeta^{x_k}}{1-\zeta^{x_k}}=\frac1{[\theta]_{x_k}^{\chi_y}}$,
where $[u]_m^{\chi_y}=\frac{(1+yu)^m-(1-u)^m}{(1+yu)^m+y(1-u)^m}$ (for
$y\neq-1$) is the $m$th power in the formal group law corresponding to the
$\chi_y$-genus (see formula~\thetag{$2_{\chi_y}$}) and
$$
\theta=\hat f_{\chi_y}\left(-\frac{2\pi i}{p}\right)=
\frac{1-e^{2\pi i/p}}{1+ye^{2\pi i/p}}=\frac{1-\zeta}{1+y\zeta}.
$$
Then $\zeta=\frac{1-\theta}{1+y\theta}$.
In what follows we assume that $y\in\Bbb Z_p$ and $y\neq-1\mod p$. The
case $y=-1\mod p$ corresponds to the Euler number which has already
been considered above. Note that for $y=0,1$ we get $\theta=1-\zeta$ and
$\theta=\frac{1-\zeta}{1+\zeta}$ respectively. This coincides with the
corresponding values for the Todd genus and the~$L$-genus obtained above.

It is easy to see that the minimal polynomial for the element
$\theta=\frac{1-\zeta}{1+y\zeta}$ is
$$
\varphi_{\theta}(u)=\frac{(1+yu)^p-(1-u)^p}{(1+y)u}.
$$
Indeed,
$$
0=\frac{\zeta^p-1}{\zeta-1}=
\frac{\left(\frac{1-\theta}{1+y\theta}\right)^p-1}
{\frac{1-\theta}{1+y\theta}-1}=
\frac1{(1+y\theta)^{p-1}}\varphi_{\theta}(\theta).
$$
Hence $\varphi_{\theta}(\theta)=0$, and $\varphi_{\theta}(u)$ is a
polynomial of degree $p-1$ with leading term $u^{p-1}$.

Further, we have
$$
\split
\prod_{k=1}^n\frac{1+y\zeta^{x_k}}{1-\zeta^{x_k}}&=
\prod_{k=1}^n\frac1{\hat f_{\chi_y}
\bigl(x_k\hat g_{\chi_y}(\theta)\bigr)}=
\prod_{k=1}^n\frac1{[\theta]^{\chi_y}_{x_k}}
\\
&=
\frac1{\theta^n}\prod_{k=1}^n\frac{\theta}{[\theta]^{\chi_y}_{x_k}}=:
\frac1{\theta^n}\sum_{i=0}^{\infty}A_i\theta^i=\frac1{\theta^n}A(\theta),
\endsplit
$$
where $A(u):=\prod_{k=1}^n\frac{u}{[u]^{\chi_y}_{x_k}}=
\sum_{i=0}^{\infty}A_iu^i$ with $A_i\in\Bbb Z_p$ being $p$-adic integers.
It follows from Lemma~5.4 that $\Tr(\theta^k)\simeq0$ for any $k>0$.

\remark{Remark}
The presentation of $\theta$ as
$\theta=C_0+C_1\zeta+\cdots+C_{p-1}\zeta^{p-1}$, \ $C_i\in\Bbb Z_p$,
which is used in the proof of Lemma~5.4, is obtained as follows. Since
$$
1+y^p=1+(y\zeta)^p=(1+y\zeta)\bigl((y\zeta)^{p-1}-(y\zeta)^{p-2}+
(y\zeta)^{p-3}-\cdots-y\zeta+1\bigr),
$$
we have
$$
\split
\theta&=\frac{1-\zeta}{1+y\zeta}=\frac{(1-\zeta)\bigl((y\zeta)^{p-1}-
(y\zeta)^{p-2}+(y\zeta)^{p-3}-\cdots-y\zeta+1\bigr)}{1+y^p}
=\frac1{1+y^p}
\\
&\qquad\times
\Bigl(\bigl(y^{p-1}+y^{p-2}\bigr)\zeta^{p-1}-
(y^{p-2}+y^{p-3})\zeta^{p-2}+\cdots-(y+1)\zeta+1-y^{p-1}\Bigr).
\endsplit
$$
Since $y^p+1\equiv y+1\neq0\mod p$, we obtain
$C_i=(-1)^i\frac{y^i+y^{i-1}}{1+y^p}\in\Bbb Z_p$ \ $(i>0)$, \
$C_0=\frac{1-y^{p-1}}{1+y^p}\in\Bbb Z_p$.
\endremark

\bigskip

So, the number-theoretical trace we are interested in is
$$
\Tr\biggl(\,\prod_{k=1}^n\frac{1+y\zeta^{x_k}}{1-\zeta^{x_k}}\biggr)
\simeq\Tr\biggl(\frac1{\theta^n}\sum_{i=0}^nA_i\theta^i\biggr)=
\sum_{i=0}^n\bigl(A_i\Tr(\theta^{i-n})\bigr)
=\bigl\langle A(u)B(u)\bigr\rangle_n,
$$
\pagebreak

\noindent
where $B_s:=\Tr\theta^{-s}$, \ $B(u):=\sum_{i=0}^{\infty}B_iu^i$. As before,
$B(u)=(p-1)+u\Tr\left(\frac1{\theta-u}\right)$, \
$\Tr\left(\frac1{\theta-u}\right)=-\frac{\varphi'_{\theta}(u)}
{\varphi_{\theta}(u)}$. Hence, the following formulae hold:
$$
\gather
{\align
\varphi_{\theta}(u)&=\frac{(1+yu)^p-(1-u)^p}{(1+y)u},
\\
\varphi'_{\theta}(u)&=\frac{py(1+yu)^{p-1}+p(1-u)^{p-1}}{(1+y)u}-
\frac{(1+yu)^p-(1-u)^p}{(1+y)u^2},
\\
\frac{\varphi'_{\theta}(u)}{\varphi_{\theta}(u)}&=
p\frac{y(1+yu)^{p-1}+(1-u)^{p-1}}{(1+yu)^p-(1-u)^p}-\frac1u,
\endalign}
\\
{\align
B(u)&=(p-1)-u\frac{\varphi'_{\theta}(u)}{\varphi_{\theta}(u)}
\\
&=
p-up\frac{y(1+yu)^{p-1}+(1-u)^{p-1}}{(1+yu)^p-(1-u)^p}=
p\frac{(1+yu)^{p-1}-(1-u)^{p-1}}{(1+yu)^p-(1-u)^p}.
\endalign}
\endgather
$$
Therefore,
$$
\split
AB_{\chi_y}(x_1,\dots,x_n)&=-\Tr\biggl(\,\prod_{k=1}^n
\frac{1+y\zeta^{x_k}}{1-\zeta^{x_k}}\biggr)
\\
&=
-\biggl\langle p\frac{(1+yu)^{p-1}-(1-u)^{p-1}}{(1+yu)^p-(1-u)^p}
\prod_{k=1}^n\frac{u}{[u]^{\chi_y}_{x_k}}\biggr\rangle_n.
\endsplit
$$
We deduce that
$$
\gather
{\align
B(u)&=p\frac{(1+yu)^{p-1}-(1-u)^{p-1}}{(1+yu)^p-(1-u)^p}=
p\frac{(1-u)(1+yu)^p-(1+yu)(1-u)^p}
{(1-u)(1+yu)\bigl((1+yu)^p-(1-u)^p\bigr)}
\\
&=
\frac p{(1-u)(1+yu)}-\frac{pu}{(1-u)(1+yu)\frac{(1+yu)^p-(1-u)^p}
{(1+yu)^p+y(1-u)^p}}
\\
&\simeq
-\frac{pu}{[u]_p^{\chi_y}}\frac1{(1-u)(1+yu)},
\endalign}
\\
\frac1{(1-u)(1+yu)}=\sum_{m=0}^{\infty}\frac{1+(-1)^my^{m+1}}{1+y}u^m,
\\
{\align
AB_{\chi_y}\left(x_1,\dots,x_n\right)&\simeq\biggl\langle
\frac{pu}{[u]_p^{\chi_y}}\frac1{(1-u)(1+yu)}
\prod_{k=1}^n\frac{u}{[u]^{\chi_y}_{x_k}}\biggr\rangle_n
\\
&=
\biggr\langle\frac{pu}{[u]_p^{\chi_y}}\prod_{k=1}^n
\frac{u}{[u]^{\chi_y}_{x_k}}\sum_{m=0}^{\infty}
\frac{1+(-1)^my^{m+1}}{1+y}u^m\biggr\rangle_n.
\endalign}
\endgather
$$
We finally obtain the following formulae for the Atiyah--Bott functions
$AB_{\chi_y}(x_1,\dots$ $\dots,x_n)=
-\Tr\left(\,\prod_{k=1}^n\frac{1+y\zeta^{x_k}}{1-\zeta^{x_k}}\right)$
of the $\chi_y$-genus ($y\neq-1\mod p$):
$$
\aligned
AB_{\chi_y}(x_1,\dots,x_n)&=-\biggl\langle
p\frac{(1+yu)^{p-1}-(1-u)^{p-1}}{(1+yu)^p-(1-u)^p}
\prod_{k=1}^n\frac{u}{[u]^{\chi_y}_{x_k}}\biggr\rangle_n,
\\
AB_{\chi_y}(x_1,\dots,x_n)&\simeq\sum_{m=0}^n\frac{1+(-1)^my^{m+1}}{1+y}
\biggl\langle\frac{pu}{[u]^{\chi_y}_p}
\prod_{k=1}^n\frac{u}{[u]^{\chi_y}_{x_k}}\biggr\rangle_{n-m}.
\endaligned
\tag{30}
$$
Here $[u]^{\chi_y}_m=\frac{(1+yu)^m-(1-u)^m}{(1+yu)^m+y(1-u)^m}$ is the~$m$th
power in the formal group law corresponding to the $\chi_y$-genus.
The $\chi_y$-genus itself is the calculated as follows:
$$
\split
\chi_y(M)&\equiv\sum_{j=1}^qAB_{\chi_y}
\bigl(x_1^{(j)},\dots,x_n^{(j)}\bigr)
\\
&\equiv\sum_{j=1}^q\sum_{m=0}^n
\frac{1+(-1)^my^{m+1}}{1+y}\biggl\langle\frac{pu}{[u]^{\chi_y}_p}
\prod_{k=1}^n\frac{u}{[u]^{\chi_y}_{x_k^{(j)}}}\biggr\rangle_{n-m}\mod p.
\endsplit
$$
Formulae \thetag{13}, \thetag{17} for the Todd
genus and the~$L$-genus are obtained from this formula by
substituting $y=0$ and $y=1$ respectively.

\head
\S\,7. The Conner--Floyd equations and the calculation of\\
Hirzebruch genera in terms of invariants of the action
\endhead

Here we consider the connection of the results in \S\,5 with the
so-called Conner--Floyd equations, which were introduced by Novikov
in~\cite{12},~\cite{13}. (Similar relations were also obtained
in~\cite{8},~\cite{11}.) Namely, it was shown there that
the sets $x_1^{(j)},\dots,x_n^{(j)}$, \ $x_k^{(j)}\in\Bbb Z/p$
are the sets of weights for some action of $\Bbb Z/p$ on a
manifold~$M^{2n}$ if and only if they satisfy the following
Conner--Floyd equations:
$$
\sum_{j=1}^q\biggl\langle\frac{pu}{[u]_p}\prod_{k=1}^n
\frac{u}{[u]_{x_k^{(j)}}}\biggr\rangle_m\simeq0,\qquad
m=0,\dots,n-1.
\tag{31}
$$
Here $[u]_m$ is the $m$th power in the universal formal group law of geometric
cobordisms (cf.~\cite{4},~\cite{12}). Applying a Hirzebruch genus
$\varphi\:\Omega_U\to\Lambda$, we obtain the Conner--Floyd equations
corresponding to~$\varphi$:
$$
\sum_{j=1}^q\biggl\langle\frac{pu}{[u]_p^{\varphi}}\prod_{k=1}^n
\frac{u}{[u]_{x_k^{(j)}}^{\varphi}}\biggr\rangle_m\simeq0,\qquad
m=0,\dots,n-1,
\tag{32}
$$
where $[u]_m^{\varphi}$ is the $m$th power in the formal group
law corresponding to the genus~$\varphi$.

The following formula for the Todd genus was deduced from
cobordism theory in~\cite{5}:
$$
\td(M)\simeq\sum_{j=1}^q\biggl\langle\frac{pu}{[u]_p^{\td}}
\prod_{k=1}^n\frac{u}{[u]_{x_k^{(j)}}^{\td}}\biggr\rangle_n.
\tag{33}
$$
Furthermore, it was shown there that formula~\thetag{33} is exactly the
difference between the formula
$$
\td(M)\simeq-\sum_{j=1}^q\Tr\biggl(\,\prod_{k=1}^n
\frac1{1-\zeta^{x_k^{(j)}}}\biggr)\simeq\sum_{j=1}^q\sum_{m=0}^n
\biggl\langle\frac{pu}{[u]_p^{\td}}
\prod_{k=1}^n\frac{u}{[u]_{x_k^{(j)}}^{\td}}\biggr\rangle_m,
$$
deduced from the Atiyah--Bott theorem (see~\S\,4.2) and the sum of the
Conner--Floyd equations~\thetag{32} for the Todd genus.
(Here $\zeta=e^{2\pi i/p}$, \ $[u]_m^{\td}=1-(1-u)^m$.)

Below we generalize this result considering the case of an arbitrary genus
$\varphi$ that satisfies the hypotheses of Theorem~5.1.

\proclaim{'heorem 7.1}
The difference between the formula of theorem~{\rm 5.1},
$$
\varphi(M)\simeq-\sum_{j=1}^q\Tr\biggl(\,\prod_{k=1}^n
\frac1{[\theta]_{x_k^{(j)}}^{\varphi}}\biggr),\qquad
\theta=f_{\varphi}\left(-\frac{2\pi i}p\right)
$$
and the sum of the Conner-Floyd equations~\thetag{32} for the genus~$\varphi$
with some~$p$-adic integer coefficients gives the following formula for the
genus~$\varphi${\rm:}
$$
\varphi(M)\simeq\sum_{j=1}^q\biggl\langle\frac{pu}{[u]_p^{\varphi}}
\prod_{k=1}^n\frac{u}{[u]_{x_k^{(j)}}^{\varphi}}\biggr\rangle_n.
$$
\endproclaim

\demo{Proof}
We have
$$
\prod_{k=1}^n\frac1{[\theta]_{x_k}^{\varphi}}=
\frac1{\theta^n}\prod_{k=1}^n\frac{\theta}{[\theta]_{x_k}^{\varphi}}=:
\frac1{\theta^n}\sum_{i=0}^{\infty}A_i\theta^i=\frac1{\theta^n}A(\theta),
$$
where $A(u):=\prod_{k=1}^n\frac u{[u]^{\varphi}_{x_k}}=
\sum_{i=0}^{\infty}A_iu^i$, \ $A_i\in\Bbb Z_p$ are~$p$-adic integers.
Therefore,
$$
\split
\Tr\biggl(\,\prod_{k=1}^n\frac1{[\theta]_{x_k}^{\varphi}}\biggr)&=
\Tr\biggl(\frac1{\theta^n}\sum_{i=0}^{\infty}A_i\theta^i\biggr)\simeq
\Tr\biggl(\frac1{\theta^n}\sum_{i=0}^{n}A_i\theta^i\biggr)
\\
&=\sum_{i=0}^{n}A_i\Tr\theta^{i-n}=
\biggl\langle\,\prod_{k=1}^n\frac u{[u]_{x_k}^{\varphi}}
B(u)\biggr\rangle_n,
\endsplit
$$
where $B_s:=\Tr\theta^{-s}$, \ $B(u):=\sum_{i=0}^{\infty}B_iu^i$.

Now we introduce a power series $h(u)$ as follows:
$$
h(u)=p\frac{[u]_p^{\varphi}-u}{B(u)[u]_p^{\varphi}}.
\tag{34}
$$
Then $h(u)$ is a series with~$p$-adic integer coefficients beginning
with~$1$. Indeed,
$$
h(0)=p\frac{[u]_p^{\varphi}-u}{B(u)[u]_p^{\varphi}}\biggr|_{u=0}=
p\frac{1-\frac u{[u]_p^{\varphi}}}{B(u)}\biggr|_{u=0}=
\frac{p\left(1-\frac1p\right)}{p-1}=1,
$$
since $B(0)=B_0=\Tr\theta^0=p-1$, \ $[u]_p^{\varphi}=pu+\cdots$. It follows
from~\thetag{34} that
$$
B(u)=\frac p{h(u)}-\frac{pu}{[u]_p^{\varphi}}\frac1{h(u)}.
$$
Hence,
$$
B(u)\simeq-\frac{pu}{[u]_p^{\varphi}}\frac1{h(u)}=
-\frac{pu}{[u]_p^{\varphi}}\biggl(1+\sum_{i=1}^{\infty}H_iu^i\biggr),
$$
where the $H_i$ are the coefficients of the series $\frac1{h(u)}$. Thus,
$$
\gather
{\align
\Tr\biggl(\,\prod_{k=1}^n\frac1{[\theta]_{x_k}^{\varphi}}\biggr)&\simeq
-\biggl\langle\frac{pu}{[u]_p^{\varphi}}\biggl(1+\sum_{i=1}^{\infty}
H_iu^i\biggr)\prod_{k=1}^n\frac u{[u]_{x_k}^{\varphi}}\biggr\rangle_n
\\
&=
-\biggl\langle\frac{pu}{[u]_p^{\varphi}}\prod_{k=1}^n
\frac u{[u]_{x_k}^{\varphi}}\biggr\rangle_n-\sum_{m=0}^{n-1}H_{n-m}
\biggl\langle\frac{pu}{[u]_p^{\varphi}}\prod_{k=1}^n
\frac u{[u]_{x_k}^{\varphi}}\biggr\rangle_m,
\endalign}
\\
{\align
\varphi(M)&\simeq-\sum_{j=1}^q\Tr\biggl(\,\prod_{k=1}^n
\frac1{[\theta]_{x_k^{(j)}}^{\varphi}}\biggr)
\\
&\simeq
\sum_{j=1}^q\biggl\langle\frac{pu}{[u]_p^{\varphi}}
\prod_{k=1}^n\frac u{[u]_{x_k^{(j)}}^{\varphi}}\biggr\rangle_n+
\sum_{m=0}^{n-1}H_{n-m}\Biggl(\,\sum_{j=1}^q\biggl\langle
\frac{pu}{[u]_p^{\varphi}}\prod_{k=1}^n
\frac u{[u]_{x_k^{(j)}}^{\varphi}}\biggr\rangle_m\Biggr).
\endalign}
\endgather
$$
This proves the theorem.
\enddemo

\example{Example 7.1}
Consider the Todd genus $\td(M)$. Then
$[u]_p^{\td}=1-(1-u)^p$, \ $B(u)=p\frac{1-(1-u)^{p-1}}{1-(1-u)^p}$
(see~\S\,4.2). Therefore,
$$
h(u)=\frac{1-(1-u)^p-u}{1-(1-u)^{p-1}}=1-u
$$
(see formula~\thetag{34}).
\endexample

\example{Example 7.2}
Consider the~$L$-genus $L(M)$. Then
$$
[u]_p^L=\frac{(1+u)^p-(1-u)^p}{(1+u)^p+(1-u)^p},\quad
B(u)=p\frac{(1+u)^{p-1}-(1-u)^{p-1}}{(1+u)^p-(1-u)^p}
$$
(see~\S\,4.3). Therefore,
$$
\split
h(u)&=\frac{\frac{(1+u)^p-(1-u)^p}{(1+u)^p+(1-u)^p}-u}
{\frac{(1+u)^{p-1}-(1-u)^{p-1}}{(1+u)^p+(1-u)^p}}
\\
&=
\frac{(1+u)^p-(1-u)^p-u(1+u)^p-u(1-u)^p}
{(1+u)^{p-1}-(1-u)^{p-1}}={(1+u)(1-u)},
\endsplit
$$
which is in accordance with the calculations from \S\,4.2, 4.3.
\endexample

\example{Example 7.3}
Consider the $\hat A$-genus $\hat A(M)$. Then
$$
[u]_p^{\hat A}=2\sh\left(p\arsh\frac u2\right),\quad
B(u)=p\frac{\sh\left((p-1)\arsh\left(\frac u2\right)\right)}
{\cosh\left(\arsh\left(\frac u2\right)\right)
\sh\left(p\arsh\left(\frac u2\right)\right)}
$$
(see~\S\,6.1). Therefore,
$$
\split
h(u)&=\frac{\left(2\sh\left(\arsh\left(\frac u2\right)\right)-u\right)
\cosh\left(\arsh\left(\frac u2\right)\right)}
{2\sh\left((p-1)\arsh\left(\frac u2\right)\right)}
\\
&=
\frac{\left(\sh\left(\arsh\left(\frac u2\right)\right)-
\sh\left(\arsh\left(\frac u2\right)\right)\right)
\cosh\left(\arsh\left(\frac u2\right)\right)}
{\sh\left((p-1)\arsh\left(\frac u2\right)\right)}
\\
&=
\frac{2\sh\left(\frac{p-1}2\arsh\left(\frac u2\right)\right)
\cosh\left(\frac{p+1}2\arsh\left(\frac u2\right)\right)
\cosh\left(\arsh\left(\frac u2\right)\right)}
{2\sh\left(\frac{p-1}2\arsh\left(\frac u2\right)\right)
\cosh\left(\frac{p-1}2\arsh\left(\frac u2\right)\right)}
\\
&=
\frac{\cosh\left(\frac{p+1}2\arsh\left(\frac u2\right)\right)
\cosh\left(\arsh\left(\frac u2\right)\right)}
{\cosh\left(\frac{p-1}2\arsh\left(\frac u2\right)\right)}.
\endsplit
$$
\endexample

\head
\S\,8. The elliptic genus for manifolds with $\Bbb Z/p$-action
\endhead

Let $M^{2n}$ be a $2n$-dimensional real orientable manifold with a complex
structure in its stable tangent bundle.

\definition{Definition 8.1 \rm(see~\cite{7})}
A Hirzebruch genus
$\varphi(M^{2n})=\left(\,\prod_{i=1}^{n}\frac{x_i}{f(x_i)}[M^{2n}]\right)$
is called the {\it elliptic genus}
if~$f$ satisfies one of the following equivalent conditions:
\medskip

1) $\dsize {f'}^2=1-2\delta f^2+\varepsilon f^4$, \ $f(0)=0$;

2) $\dsize f(u+v)=\frac{f(u)f'(v)+f'(u)f(v)}{1-\varepsilon f(u)^2f(v)^2}$.
\enddefinition

Let us consider the lattice $L=2\pi i(\Bbb Z\tau+\Bbb Z)$ in~$\Bbb C$, with
$\Ima\tau>0$, and put $L'=L\backslash\{0\}$. The Weierstra\ss\
$\wp$-function is
$\wp(z)=\frac1{z^2}+\sum_{\omega\in L'}\left(\frac1{(z-\omega)^2}-
\frac1{\omega^2}\right)$. It satisfies the following differential equation:
$$
\wp'(z)^2=4\wp(z)^3-g_2\wp(z)-g_3=
4\bigl(\wp(z)-e_1\bigr)\bigl(\wp(z)-e_2\bigr)
\bigl(\wp(z)-e_3\bigr),
$$
where $e_1=\wp(\pi i)$, \ $e_2=\wp(\pi i\tau)$, \
$e_3=\wp\bigl(\pi i(\tau+1)\bigr)$ are the zeros of the derivative $\wp'$.
The function $f(z)=1\big/\!\sqrt{\wp(z)-e_1}$ (that is, $f(z)=\sn(z)$, the
elliptic sine) satisfies the conditions of Definition~8.1 for
$\delta=-\frac32e_1$, \ $\varepsilon=(e_1-e_2)(e_1-e_3)$ (see~\cite{7}).
Thus it gives rise to the elliptic genus. However,
$\varepsilon=(e_1-e_2)(e_1-e_3)\neq0$ and
$\delta^2-\varepsilon=\frac14(e_2-e_3)^2\neq0$,
since $e_1$, $e_2$,~$e_3$ are all different. The degenerate cases
$e_1=e_2$ and $e_1=e_3$ ($\varepsilon=0$)
correspond to the $\hat A$-genus.
(Putting $\delta=-1/8$, we obtain
${f'}^2=1+\frac14f^2$, that is, $f(x)=2\sh(x/2)$.)
The degenerate case $e_2=e_3$ ($\delta^2-\varepsilon=0$) corresponds to
the $L$-genus. (Putting
$\delta=\varepsilon=1$, we obtain ${f'}^2=(1-f^2)^2$, that is, $f(x)=\th x$.)

The differential equation defining the function $f(x)$ for the
elliptic genus implies that if we put $\deg\varepsilon=4$ and
$\deg\delta=2$, then the coefficient of~$x^{2k}$ in the series $f(x)$ becomes
a weighted homogeneous polynomial of degree~$2k$ in~$\delta$
and~$\varepsilon$ with coefficients in the ring $\Bbb Z\left[\frac12\right]$.
Therefore, the elliptic genus
$\varphi(M^{2n})=\left(\,\prod_{i=1}^n\frac{x_i}{f(x_i)}\right)[M^{2n}]$
is a homogeneous polynomial of degree~$2n$ in~$\delta$
and~$\varepsilon$. So $\varphi$ can be regarded as a homomorphism from
the complex cobordism ring~$\Omega_U$ to the ring
$\Bbb Z\left[\frac12\right][\delta,\varepsilon]$ (or to the ring
$\Bbb Z_p[\delta,\varepsilon]$, \ $p>2$).

The function $f(\tau,x)$ is an elliptic function for the
sublattice $\tilde L=2\pi i(\Bbb Z\cdot 2\tau+\Bbb Z)\subset L$ of index 2.
The divisor of this function is
$(0)+(\pi i\cdot 2\tau)-(\pi i)-\bigl(\pi i(1+2\tau)\bigr)$. We have the
following decomposition of $f(\tau,x)$ into an infinite product
(see~\cite{7}):
$$
f(\tau,x)=2\frac{1-e^{-x}}{1+e^{-x}}\prod_{k=1}^{\infty}
\frac{(1-q^ke^x)(1-q^ke^{-x})(1+q^k)^2}
{(1+q^ke^x)(1+q^ke^{-x})(1-q^k)^2},\quad q=e^{2\pi i\tau}.
\tag{35}
$$

{\tolerance=2000
Now we consider the following object:
$$
\Biggl\{L_i^{(q)}=\sum_{p=0}^n\Lambda^{p,i}\otimes\biggl(\,
\bigotimes_{k=1}^{\infty}S_{q^k}T_{\Bbb C}M\otimes
\bigotimes_{k=1}^{\infty}\Lambda_{q^k}T_{\Bbb C}M\biggr)\Biggr\},
\tag{36}
$$
where $T_{\Bbb C}M$~is the complexification of the tangent bundle to $M$,
$\Lambda_tE:=$ $\sum_{k=0}^{\infty}(\Lambda^kE)t^k$, \
$S_tE:=\sum_{k=0}^{\infty}(S^kE)t^k$. Then $L^{(q)}$ is a power series
in~$q$ whose coefficients
are elliptic complexes associated to $TM$. Its
index is a well-known power series in~$q$, namely, the so-called twisted
signature (see~\cite{7}):
$\inda(L^{(q)})=\sign\left(M,\bigotimes_{k=1}^{\infty}
S_{q^k}T_{\Bbb C}M\otimes\bigotimes_{k=1}^{\infty}
\Lambda_{q^k}T_{\Bbb C}M\right)$. It follows from the Atiyah--Singer
theorem~2.1 and Lemma~2.3 that this index is
$$
\inda(L^{(q)})=\Biggl(\,\prod_{j=1}^n\biggl(x_j
\frac{1+e^{-x_j}}{1-e^{-x_j}}\prod_{k=1}^{\infty}
\frac{(1+q^ke^{x_j})(1+q^ke^{-x_j})}{(1-q^ke^{x_j})(1-q^ke^{-x_j})}
\biggr)\Biggr)[M^{2n}].
\tag{37}
$$
It is clear from this formula that $\inda(L^{(q)})$ is a power series in~$q$
with integer coefficients and with constant term
(the coefficient of~$q^0$) \ $L(M)=\sign M$.
Comparing this expression with~\thetag{35} and taking into account that
$$
\frac12\prod_{k=1}^{\infty}\frac{(1-q^k)^2}{(1+q^k)^2}=\varepsilon^{1/4}
\tag{38}
$$
(see~\cite{7}), we obtain
$$
\varphi(M)=\sign\biggl(M,\bigotimes_{k=1}^{\infty}S_{q^k}
T_{\Bbb C}M\otimes\bigotimes_{k=1}^{\infty}\Lambda_{q^k}
T_{\Bbb C}M\biggr)\varepsilon^{n/4}.
\tag{39}
$$

}
Suppose that an operator $g$, $g^p=\nomathbreak1$, acts on a
manifold~$M^{2n}$ with finitely many fixed points
$\Cal P_1,\dots,\Cal P_r$. According to the recipe from~\S\,5,
the equivariant index $\inda\bigl(g,L^{(q)}\bigr)$ of the complex~$L^{(q)}$
equals the sum of the contribution functions $\sigma_{L^{(q)}}(\Cal P_j)$ of
the fixed points. These contributions are obtained from formula~\thetag{37} by
replacing~$M^{2n}$ by~$\Cal P_j$,
\ $x_1x_2\dots x_n$ by~1 and~$x_k$ by $-\frac{2\pi i}px_k^{(j)}$:
$$
\sigma_{L^{(q)}}(\Cal P_j)=\prod_{k=1}^n\left(
\frac{1+\zeta^{x_k^{(j)}}}{1-\zeta^{x_k^{(j)}}}\prod_{i=1}^{\infty}
\frac{\Bigl(1+q^i\zeta^{-x_k^{(j)}}\Bigr)\Bigl(1+q^i\zeta^{x_k^{(j)}}\Bigr)}
{\Bigl(1-q^i\zeta^{-x_k^{(j)}}\Bigr)\Bigl(1-q^i\zeta^{x_k^{(j)}}\Bigr)}
\right),\quad
\zeta=e^{\frac{2\pi i}p}.
$$

\pagebreak

We note that
$\frac1p\sum_{l\in\Bbb Z/p}\inda\bigl(g^l,L^{(q)}\bigr)=s(q)$ is a power
series in~$q$ whose coefficient at~$q^k$ is the alternating sum of the
dimensions of the invariant subspaces for the action of~$g$
on the cohomology of a certain complex. Namely, this complex is the
coefficient of~$q^k$ in the series~$L^{(q)}$:
$$
L^{(q)}=L+2L\otimes T_{\Bbb C}M q+L\otimes
(2T_{\Bbb C}M+T_{\Bbb C}M\otimes T_{\Bbb C}M+S^2T_{\Bbb C}M+
\Lambda^2T_{\Bbb C}M)q^2+\cdots,
$$
where $L=\sum_{p=0}^n\Lambda^{p,i}$. Therefore, $s(q)$ is a series with
integer coefficients. Furthermore, $\inda\bigl(1,L^{(q)}\bigr)=\sign\left(M,
\bigotimes_{k=1}^{\infty}S_{q^k}T_{\Bbb C}M\otimes
\bigotimes_{k=1}^{\infty}\Lambda_{q^k}T_{\Bbb C}M\right)$. Hence,
$$
\split
&\sign\biggl(M,\bigotimes_{k=1}^{\infty}S_{q^k}T_{\Bbb C}M\otimes
\bigotimes_{k=1}^{\infty}\Lambda_{q^k}T_{\Bbb C}M\biggr)
\\
&\qquad=
\inda(1,L^{(q)})=-\sum_{l=1}^{p-1}\inda(g^l,L^{(q)})+ps(q)
\\
&\qquad=
-\sum_{j=1}^r\sum_{l=1}^{p-1}\prod_{k=1}^n\left(
\frac{1+\zeta^{lx_k^{(j)}}}{1-\zeta^{lx_k^{(j)}}}\prod_{i=1}^{\infty}
\frac{\Bigl(1+q^i\zeta^{-lx_k^{(j)}}\Bigr)\Bigl(1+q^i\zeta^{lx_k^{(j)}}\Bigr)}
{\Bigl(1-q^i\zeta^{-lx_k^{(j)}}\Bigr)
\Bigl(1-q^i\zeta^{lx_k^{(j)}}\Bigr)}\right)+ps(q).
\endsplit
$$
The left-hand side of this relation belongs to the ring
$\Bbb Z\bigl[[q]\bigr]$ of power series with integer coefficients, while its
right-hand side {\it a priori}
belongs to an extension of $\Bbb Z\bigl[[q]\bigr]$, namely, to
the ring $\Bbb Z\bigl[[q]\bigr](\zeta)$, \ $\zeta^p=1$. We embed
both rings in the~$p$-adic extensions of the corresponding fields and
consider the number-theoretical trace
$\Tr\:\Bbb Q_p\{q\}(\zeta)\to\Bbb Q_p\{q\}$, where
$\Bbb Q_p\{q\}=\Bbb Q_p[[q]][q^{-1}]$ is the
field of Laurent series in~$q$ with rational~$p$-adic coefficients and
$\Bbb Q_p\{q\}(\zeta)$ is the algebraic extension of $\Bbb Q_p\{q\}$ by the
element~$\zeta$, \ $\zeta^p=1$. Then
$$
\split
&\sum_{l=1}^{p-1}\prod_{k=1}^n\left(
\frac{1+\zeta^{lx_k^{(j)}}}{1-\zeta^{lx_k^{(j)}}}\prod_{i=1}^{\infty}
\frac{\Bigl(1+q^i\zeta^{-lx_k^{(j)}}\Bigr)\Bigl(1+q^i\zeta^{lx_k^{(j)}}\Bigr)}
{\Bigl(1-q^i\zeta^{-lx_k^{(j)}}\Bigr)\Bigl(1-q^i\zeta^{lx_k^{(j)}}\Bigr)}
\right)
\\
&\qquad=
\Tr\left(\prod_{k=1}^n\left(\frac{1+\zeta^{x_k^{(j)}}}
{1-\zeta^{x_k^{(j)}}}\prod_{i=1}^{\infty}
\frac{\Bigl(1+q^i\zeta^{-x_k^{(j)}}\Bigr)
\Bigl(1+q^i\zeta^{x_k^{(j)}}\Bigr)}
{\Bigl(1-q^i\zeta^{-x_k^{(j)}}\Bigr)
\Bigl(1-q^i\zeta^{x_k^{(j)}}\Bigr)}\right)\right).
\endsplit
$$
Thus the following proposition holds:

\proclaim{Proposition 8.2}
We have the following formula for the index of the complex
$L^{(q)}$ in~\thetag{36}:
$$
\split
&\inda(L^{(q)})=\sign\biggl(M,\bigotimes_{k=1}^{\infty}S_{q^k}T_{\Bbb
C}M\otimes \bigotimes_{k=1}^{\infty}\Lambda_{q^k}T_{\Bbb C}M\biggr) \\ &\,\,
\equiv
-\sum_{j=1}^r\Tr\!\left(\prod_{k=1}^n\!\left(\!
\frac{1+\zeta^{x_k^{(j)}}}{1-\zeta^{x_k^{(j)}}}
\prod_{i=1}^{\infty}
\frac{\Bigl(1+q^i\zeta^{-x_k^{(j)}}\Bigr)\Bigl(1+q^i\zeta^{x_k^{(j)}}\Bigr)}
{\Bigl(1-q^i\zeta^{-x_k^{(j)}}\Bigr)\Bigl(1-q^i\zeta^{x_k^{(j)}}\Bigr)}
\right)\!\right)
\bmod\, p\Bbb Z\bigl[[q]\bigr].\qquad
\endsplit
\tag{40}
$$
\endproclaim

{\tolerance=2000
Using relations \thetag{38} and \thetag{39}, we thus
obtain the following formula for the elliptic genus~$\varphi(M)$:
$$
\split
\varphi(M)&\equiv-\sum_{j=1}^r\Tr\left(\,\prod_{k=1}^n\left(
\frac12\cdot\frac{1+\zeta^{x_k^{(j)}}}{1-\zeta^{x_k^{(j)}}}
\vphantom{\frac{\Bigl(1+\zeta^{x_k^{(j)}}\Bigr)}
{\Bigl(1-\zeta^{x_k^{(j)}}\Bigr)}}
\right.\right.
\\
&\qquad\left.\left.
\times
\prod_{i=1}^{\infty}\frac{\bigl(1+q^i\zeta^{-x_k^{(j)}}\bigr)
\Bigl(1+q^i\zeta^{x_k^{(j)}}\Bigr)(1-q^i)^2}
{\Bigl(1-q^i\zeta^{-x_k^{(j)}}\Bigr)
\Bigl(1-q^i\zeta^{x_k^{(j)}}\Bigr)(1+q^i)^2}\right)\right)
\mod p\Bbb Z_p\bigl[[q]\bigr].
\endsplit
$$
Taking into account the decomposition \thetag{35} of $f(\tau,x)=\sn x$, we
rewrite this formula as
$$
\varphi(M)\equiv-\sum_{j=1}^r\Tr\left(\,\prod_{k=1}^n
\frac1{f\Bigl(\tau,-\frac{2\pi ix_k^{(j)}}p\Bigr)}
\right)\mod p\Bbb Z_p\bigl[[q]\bigr],\quad q=e^{2\pi i\tau}.
$$
However, $\varphi(M)$ is a homogeneous polynomial in
$\delta$, $\varepsilon$ with coefficients in the ring
$\Bbb Z\bigl[\frac12\bigr]\subset\Bbb Z_p$. So we want to obtain the
expression for $\varphi(M)$ in the ring $\Bbb Z_p[\delta,\varepsilon]$,
not in the ring $\Bbb Z_p\bigl[[q]\bigr]$.
($\Bbb Z_p[\delta,\varepsilon]$ is the subring of $\Bbb Z_p\bigl[[q]\bigr]$.)
As before, we put $\theta=f\left(\tau,-\frac{2\pi i}p\right)$. Then
$f\left(\tau,-\frac{2\pi ix_k}p\right)=[\theta]_{x_k}$, where
$[\theta]_m=f\bigl(mf^{-1}(\theta)\bigr)$ is the $m$th power in the formal
group law corresponding to the elliptic genus. We note that
$\theta$ is an algebraic element over the field
$\Bbb Q(\delta,\varepsilon)$. Indeed,
$$
[\theta]_p=f\left(\tau,-\frac{2\pi ip}p\right)=f(\tau,-2\pi i)=0,
$$
and $[u]_p$ is a rational function of $u$ when $p$ is odd.
(For example, $[u]_3=
u\frac{3-8\delta u^2+6\varepsilon u^4-\varepsilon^2u^8}
{1-6\varepsilon u^4+8\delta\varepsilon u^6-3\varepsilon^2u^8}$.
This follows from the differential equation and the addition theorem for
$f(\tau,x)=\sn x$.) The numerator of $[u]_p$ is a polynomial in~$u$ with
coefficients $\delta$,~$\varepsilon$, and~$\theta$ is a zero of this
polynomial. Obviously, this polynomial is of degree~$p^2$ in~$u$, and its
zeros are the values of~$f$ at all~$p$-division points of the lattice
$2\pi i(\Bbb Z\cdot 2\tau+\Bbb Z)$, that is, they equal
$f\left(\tau,2\pi i\frac{2k\tau+l}p\right)$, $k,l\in\Bbb Z$. In particular,
$\theta,[\theta]_2,\dots,[\theta]_{p-1},[\theta]_p=0$ are among zeros.

}
Finally, we have
$$
\varphi(M)\equiv-\sum_{j=1}^r\Tr\biggl(\,\prod_{k=1}^n
\frac1{[\theta]_{x_k}}\biggr)\mod p\Bbb Z_p\bigl[[q]\bigr],\quad
\theta=f\left(\tau,-\frac{2\pi i}p\right),
\tag{41}
$$
where both right- and left-hand sides are polynomials in
$\delta$,~$\varepsilon$.

\proclaim{Lemma 8.3}
Let $P(\delta,\varepsilon)$ be a homogeneous polynomial in
$\delta$,~$\varepsilon$ such that
$P(\delta,\varepsilon)\equiv0\mod p\Bbb Z_p\bigl[[q]\bigr]$. Then
$P(\delta,\varepsilon)\equiv0\mod p\Bbb Z_p[\delta,\varepsilon]$, that is,
all coefficients of $P(\delta,\varepsilon)$ belong to $p\Bbb Z_p$.
\endproclaim

\demo{Proof}
The functions $\delta$,~$\varepsilon$ are modular forms of
weights~2 and~4 respectively on the subgroup
$\Gamma_0(2)\subset\SL_2(\Bbb Z)$, \
$\Gamma_0(2):=\biggl\{A\in\SL_2(\Bbb Z)\biggm|A\equiv
\pmatrix * & * \\ 0 & *\endpmatrix\pmod2\biggr\}$. They have
the following Fourier expansions near infinity (that is, for small
$q=e^{2\pi i\tau}$):
$$
\align
\delta&=-\frac32e_1=\frac14+\sum_{n=1}^{\infty}
\sum_{\Sp d\mid n \\ d\equiv1\mod2\endSp}dq^n
\\
&=\frac14+6(q+q^2+4q^3+q^4+6q^5+4q^6+\cdots),
\\
\varepsilon&=(e_1-e_2)(e_1-e_3)
\\
&=\frac1{16}+\sum_{n=1}^{\infty}\sum_{d\mid n}(-1)^dd^3q^n=
\frac1{16}-q+7q^2-28q^3+71q^4-\cdots
\endalign
$$
Furthermore, $\delta(0)=-\frac18$, \ $\varepsilon(0)=0$,  \
$\delta\left(\frac{1+i}2\right)=0$, \
$\varepsilon\left(\frac{1+i}2\right)\neq0$ (see~\cite{7}).

The proof is by induction on the degree $\deg P=2n$
of the polynomial~$P$. If $\deg P(\delta,\varepsilon)=2$, then
$P(\delta,\varepsilon)=b\delta$. Hence, $b\in p\Bbb Z_p$, because
$\delta\notin p\Bbb Z_p$. Now assume that the lemma holds for all
polynomials of degree $<2k$ and let $\deg P(\delta,\varepsilon)=2k$.
Write $P(\delta,\varepsilon)=
a\varepsilon^{k/2}+Q(\delta,\varepsilon)\varepsilon\delta+b\delta^k$,
where $\deg Q(\delta,\varepsilon)<2k$, and therefore
$Q(\delta,\varepsilon)\equiv0\mod p\Bbb Z_p[\delta,\varepsilon]$.
Evaluation of $P(\delta,\varepsilon)$ at the point $\tau=\frac{1+i}2$ shows
that $a\varepsilon\left(\frac{1+i}2\right)^{k/2}$ is divisible by $p$, whence
$a$ is divisible by $p$. Evaluation of $P(\delta,\varepsilon)$ at $\tau=0$
shows that $\left(-\frac18\right)^kb$ is divisible by $p$, whence
$b$ is divisible by $p$. Thus all the coefficients of
$P(\delta,\varepsilon)$ belong to $p\Bbb Z_p$. The lemma is proved.
\enddemo

From this lemma and \thetag{41} we obtain the following
theorem.

\proclaim{Theorem 8.4}
Let $g$ be a transversal endomorphism  acting on~$M^{2n}$
such that $g^p=1$. Then we have the following formula for the elliptic genus
$\varphi(M^{2n})${\rm:}
$$
\varphi(M)\equiv-\sum_{j=1}^r\Tr\left(\,\prod_{k=1}^n
\frac1{f\Bigl(\tau,-\frac{2\pi ix_k^{(j)}}p\Bigr)}\right)
\mod p\Bbb Z_p[\delta,\varepsilon],
\tag{42}
$$
where $\Tr$ in the right-hand side denotes the number-theoretical trace
$\Bbb Q_p(\delta,\varepsilon)(\zeta)\to\Bbb Q_p(\delta,\varepsilon)$.
\endproclaim

All constructions in~\S\,5 and~\S\,7 (in particular, Lemma~5.4 and
Theorem~7.1) can be repeated for the elliptic genus
(we just replace the ring~$\Bbb Z_p$ by
$\Bbb Z_p[\delta,\varepsilon]$ in all formulae). Thus, the difference
between formula~\thetag{42} and a certain weighted sum of the following
Conner--Floyd equations for the elliptic genus:
$$
\sum_{j=1}^r\biggl\langle\frac{pu}{[u]_p}
\prod_{k=1}^n\frac u{[u]_{x_k^{(j)}}}\biggr\rangle_m\equiv0
\mod p\Bbb Z_p[\delta,\varepsilon],\quad m=0,\dots,n-1,
$$
gives the following formula for the elliptic genus:
$$
\varphi(M^{2n})\equiv\sum_{j=1}^r\biggl\langle\frac{pu}{[u]_p}
\prod_{k=1}^n\frac u{[u]_{x_k^{(j)}}}\biggr\rangle_n
\mod p\Bbb Z_p[\delta,\varepsilon].
\tag{43}
$$
In the last two formulae, $[u]_m$ denotes the $m$th power in the formal
group law corresponding to the elliptic genus.

\medskip

To give an example, we apply the formulae for the elliptic genus to a
particular action of the group $\Bbb Z/p$. Namely, we consider the
action of $\Bbb Z/p$ on $\Bbb CP^n$ such that the generator~$g$ acts in
homogeneous coordinates as
$(z_0:z_1:\cdots:z_n)\to(\lambda_0z_0:\lambda_1z_1:\cdots:\lambda_nz_n)$, \
$\lambda_i=e^{\frac{2\pi i}py_i}$. If $y_0,\dots,y_n$ are distinct as
residues $\bmod\, p$, which implies that ${n<p}$, then this action has only
finitely many fixed points on $\Bbb CP^n$. Namely, the fixed points are $\Cal
P_j=\Bigl(0,\dots,\overset{j}\to{1},\dots,0\Bigr)$, \ $j=0,\dots,n$. In
the local coordinates $\left(\frac{z_0}{z_j},\frac{z_1}{z_j},\dots,
\widehat{\frac{z_j}{z_j}},\dots,\frac{z_n}{z_j}\right)$
near~$\Cal P_j$, the operator~$g$ acts linearly:
$\frac{z_i}{z_j}\to\frac{\lambda_iz_i}{\lambda_jz_j}=
\exp\left(\frac{2\pi i}{p}(y_i-y_j)\right)\frac{z_i}{z_j}$. Therefore,
the eigenvalues of the map~$g$ at the point~$\Cal P_j$ are
$\exp\left(\frac{2\pi i}{p}(y_i-y_j)\right)$, \ $i\neq j$, and
the corresponding weights are $x_i^{(j)}=y_i-y_j$.

Now we consider the elliptic genus~$\varphi$. It is well known that
$$
\varphi(\Bbb CP^n)=\left\langle\frac1
{\sqrt{1-2\delta u^2+\varepsilon u^4}}\right\rangle_n.
\tag{44}
$$

Indeed, it follows from formula \thetag{1}, that
$\varphi(\Bbb CP^n)=\bigl\langle g^{\prime}_{\varphi}(u)\bigr\rangle_n$,
where $g_{\varphi}(u)=f_{\varphi}^{-1}(u)$ is the logarithm of the
corresponding formal group law. But
$g^{\prime}_{\varphi}(u)=\frac1{\sqrt{1-2\delta u^2+\varepsilon u^4}}$,
since $(f^{\prime})^2=1-2\delta f^2+\varepsilon f^4$.

At the same time, it follows from~\thetag{43} that
$$
\varphi(\Bbb CP^n)\equiv\sum_{j=0}^n\biggl\langle\frac{pu}{[u]_p}
\prod_{i\neq j}\frac u{[u]_{y_i-y_j}}\biggr\rangle_n
\mod p\Bbb Z_p[\delta,\varepsilon].
$$
Hence for any set $y_0,y_1,\dots,y_n$ of distinct residues modulo $p$, where
$n<p$, we have
$$
\sum_{j=0}^n\biggl\langle\frac{pu}{[u]_p}
\prod_{i\neq j}\frac u{[u]_{y_i-y_j}}\biggr\rangle_n\equiv
\left\langle\frac1{\sqrt{1-2\delta u^2+\varepsilon u^4}}\right\rangle_n
\mod p\Bbb Z_p[\delta,\varepsilon]
$$
The function ${(1-2\delta u+u^2)}^{-1/2}$ is the generating function for
the Legendre polynomials $P_n(\delta)$, that is,
$\frac1{\sqrt{1-2\delta u+u^2}}=1+\sum_{n>0}P_n(\delta)u^n$. Hence,
$$
\frac1{\sqrt{1-2\delta u^2+\varepsilon u^4}}=
1+\sum_{n>0}P_n\left(\frac{\delta}k\right)(ku^2)^n,
$$
with $k^2=\varepsilon$. Therefore,
$$
\aligned
\left\langle\frac1{\sqrt{1-2\delta u^2+\varepsilon u^4}}
\right\rangle_{2m}&=k^mP_m\left(\frac{\delta}k\right),
\\
\sum_{j=0}^{2m}\biggl\langle\frac{pu}{[u]_p}
\prod_{i\neq j}\frac u{[u]_{y_i-y_j}}\biggr\rangle_{2m}&\equiv
k^mP_m\left(\frac{\delta}k\right)
\mod p\Bbb Z_p[\delta,\varepsilon].
\endaligned
\tag{45}
$$
In the simplest case $2m=n=p-1$, the set
$y_0-y_j,y_1-y_j,\dots,\widehat{y_j-y_j},\dots,y_{p-1}-y_j$ forms the
complete set $1,2,\dots,p-1$ of non-zero residues modulo $p$ (since the~$y_i$
are distinct). Hence, all the summands in the left-hand side of the last
formula\linebreak

\pagebreak

\noindent are equal, and we obtain
$$
\split
\left\langle\frac{p^2u^p}{u[u]_2[u]_3\cdots[u]_{p-1}[u]_p}
\right\rangle_{p-1}
&\equiv
\left\langle
\frac1{\sqrt{1-2\delta u^2+\varepsilon u^4}}\right\rangle_{p-1}
\\
&=
k^{\frac{p-1}2}P_{\frac{p-1}2}\left(\frac{\delta}k\right)
\mod p\Bbb Z_p[\delta,\varepsilon].
\endsplit
\tag{46}
$$
This formula could be also deduced from the well-known relation
$$
[u]_p\equiv P_{(p-1)/2}(\delta)u^p+\cdots\mod p\Bbb Z_p[\delta]
$$
(see~\cite{10}). Here we put $k=1$ for simplicity, and the dots stand for
the higher order terms. Indeed, we could rewrite the above relation as
$$
[u]_p=pu(1+b_1u+\cdots+b_{p-1}u^{p-1})+P_{(p-1)/2}(\delta)u^p+\cdots,
$$
with $b_1,\dots,b_{p-1}\in\Bbb Z[\delta]$. Hence,
$$
\split
&\left\langle\frac{p^2u^p}{u[u]_2\cdots[u]_{p-1}[u]_p}
\right\rangle_{p-1}
\\
&\qquad
=\biggl\langle\frac{p^2}{u(2u+\cdots)
\dots\bigl((p-1)u+\cdots\bigr)}
\\
&\qquad\qquad\times
\frac{u^p}{\bigl(pu(1+b_1u+\cdots+b_{p-1}u^{p-1})+
P_{(p-1)/2}(\delta)u^p+\cdots\bigr)}\biggr\rangle_{p-1}
\\
&\qquad
=\left\langle\frac1{(p-1)!}\frac{p^2}{p(1+c_1u+\cdots+c_{p-1}u^{p-1})
+P_{(p-1)/2}(\delta)u^{p-1}+\cdots}\right\rangle_{p-1}
\\
&\qquad
=\frac1{(p-1)!}\left\langle\frac{p}{1+c_1u+\cdots+c_{p-2}u^{p-2}+
\bigl(c_{p-1}+P_{(p-1)/2}(\delta)/p\bigr)u^{p-1}}\right\rangle_{p-1},
\endsplit
$$
where $c_1,\dots,c_{p-1}$ belong to $\Bbb Z[\delta]$. Calculating the
latter expression modulo $p$, we get
$$
\left\langle\frac{p^2u^p}{u[u]_2[u]_3\cdots[u]_{p-1}[u]_p}
\right\rangle_{p-1}\equiv-\frac1{(p-1)!}P_{(p-1)/2}(\delta)\equiv
P_{(p-1)/2}(\delta)\mod p\Bbb Z_p[\delta],
$$
since $(p-1)!\equiv-1\mod p$, as required.

\head
\S\,9. Generalization to $\Bbb Z/p$-actions having\\
fixed submanifolds with trivial normal bundle
\endhead

Suppose that an operator~$g$, $g^p=1$, acts on a stably complex
manifold $M^{2n}$. Let the fixed point set~$M^g$ be written as the union
of connected fixed submanifolds: $M^g=\bigcup_{\nu}M_{\nu}^g$. Also
suppose that all the~$M_{\nu}^g$ have trivial normal bundle in~$M$.

Let $\varphi$ be a Hirzebruch genus that can be calculated as the
index of an elliptic complex $E_\varphi$ of bundles associated to the
tangent bundle~$TM$. To avoid misunderstanding, we shall
denote the first Chern classes of the ``virtual" line subbundles of $TM$
by~$z_i$. Then $c(TM)= (1+z_1)\dots(1+z_n)$.

Let $Y=M^g_k$ be one of the fixed submanifolds for the action of~$g$.
According to formulae~\thetag{20} from~\S\,5, we have:
$$
\gather
TM|_Y=\bigoplus_{j}N_{\lambda_j},\qquad
c(N_{\lambda_j})=\prod_{i=1}^{d_{\lambda_j}}\Bigl(1+z_i^{(\lambda_j)}\Bigr),
\\
c(TM_Y)=\prod_{j}\prod_{i=1}^{d_{\lambda_j}}\Bigl(1+z_i^{(\lambda_j)}\Bigr)=
\prod_{i=1}^n(1+z_i),
\endgather
$$
where $N_{\lambda_j}$ is the subbundle of $TM|_Y$ corresponding to the
eigenvalue~$\lambda_j$ of the differential of~$g$. Here $\lambda_j^p=1$
(hence, $\lambda_j=e^{2\pi ix_j/p}$) and $d_{\lambda_j}=\dim N_{\lambda_j}$.
Obviously, $N_1=TY$ is the tangent bundle to~$Y$.
According to the results of Atiyah and Singer~\cite{3} described in~\S\,5, the
equivariant index $\inda(g,E_{\varphi})$ of the complex~$E_{\varphi}$ can be
computed as $\inda(g,E_{\varphi})=\sum_{\nu}\sigma(M^g_{\nu})$. We also have
the following ``recipe" for calculating the fixed point contribution
functions $\sigma(Y)$. In the formula~\thetag{21} for the index
of~$E_{\varphi}$, replace~$M$ by~$Y$ and $e^{z_i^{(\lambda_j)}}$
by $\lambda_j^{-1}e^{z_i^{(\lambda_j)}}$. Equivalently, since
$\lambda_j=e^{2\pi ix_j/p}$, one should replace $z_i^{(\lambda_j)}$ by
$z_i^{(\lambda_j)}-2\pi ix_j/p$.

Since $\inda(E_{\varphi})=\varphi(M)$, we have
$$
\inda(E_{\varphi})=\biggl(\,\prod_{i=1}^n\frac{z_i}{f(z_i)}\biggr)[M]=
\Biggl(\,\prod_j\prod_{i=1}^{d_{\lambda_j}}
\frac{z_i^{(\lambda_j)}}{f\bigl(z_i^{(\lambda_j)}\bigr)}\Biggr)[M].
$$
Hence, our ``recipe" shows that
$$
\sigma(Y)=\Biggl(\,\prod_{\lambda_j\neq1}\prod_{i=1}^{d_{\lambda_j}}
\frac1{f\bigl(z_i^{(\lambda_j)}-2\pi ix_j/p\bigr)}
\prod_{i=1}^{d_1}\frac{z_i^{(1)}}{f\bigl(z_i^{(1)}\bigr)}\Biggr)[Y],
$$
because $c_n(M)=\prod_j\prod_{i=1}^{d_{\lambda_j}}z_i^{(\lambda_j)}$
reduces to $c_n(Y)=\prod_{i=1}^{d_1}z_i^{(1)}$ and all the weights~$x_k$
corresponding to the eigenvalue $\lambda=1$ are zero. Furthermore,
since~$Y$ has trivial normal bundle in~$M$,
we have $z_i^{(\lambda_j)}=0$ for $\lambda_j\neq1$. Therefore,
$\sigma(Y)=\prod_j\frac1{f(-2\pi ix_j/p)}\varphi(Y)$, where the~$x_j$ are the
weights corresponding to the eigenvalues
$\lambda_j=e^{2\pi ix_j/p}\neq1$, that is, the $x_j$ are non-zero modulo $p$.
Finally, we have the following formula for the equivariant index:
$$
\inda(g,E_{\varphi})=\sum_{\nu}\Biggl(\,\prod_j
\frac1{f\bigl(-2\pi ix_j^{(\nu)}\bigr)}\varphi(M^g_{\nu})\Biggr).
$$

Theorem~5.1 and Proposition~5.3 can naturally be extended to the case of
actions having fixed submanifolds with trivial normal bundle. We get
the following statement.

\proclaim{Theorem 9.1}
Let $\varphi$ be a Hirzebruch genussuch that there is an
elliptic complex of bundles associated to $TM$ whose index is equal to
$\varphi(M)$. Then we have the\linebreak

\pagebreak

\noindent following formula for $\varphi(M)${\rm:}
$$
\split
\varphi(M)&\equiv-\sum_{\nu}\Tr\Biggl(\,\prod_j
\frac1{f\bigl(-2\pi ix_j^{(\nu)}/p\bigr)}\Biggr)\varphi(M_{\nu}^g)
\\
&=
-\sum_{\nu}\Tr\biggl(\,\prod_j
\frac1{[\theta]^{\varphi}_{x_j^{\nu}}}\biggr)
\varphi(M_{\nu}^g)\mod p,
\endsplit
$$
where $M^g=\bigcup_{\nu}M^g_{\nu}$ is the fixed point set,
$\theta=f_{\varphi}(-2\pi i/p)$, \ $[\theta]_k$ is the $k$-th power in the
corresponding formal group, and
$\Tr\:\Bbb Q(e^{2\pi i/p})\to\Bbb Q$ is the number-theoretical trace.
{\rm(}If the genus~$\varphi$ takes its values in some ring~$\Lambda$
instead of~$\Bbb Z$, replace~$\Bbb Q$ by the
corresponding quotient field\/{\rm)}.
\endproclaim

The author is grateful to Prof. V.\,M.~Buchstaber and Prof.
S.\,P.~Novikov for suggesting the problem and wishes to express
special thanks to V.\,M.~Buchstaber for useful recommendations and
stimulating discussions.


\tolerance=2000

\enddocument